\documentclass{article}
\usepackage{amsmath, amsthm, amssymb}
\usepackage[margin=1in]{geometry}
\usepackage[labelformat=empty]{caption}

\newif\ifpdf \pdffalse \ifx\pdfoutput\undefined\else\ifx\pdfoutput\relax\else\ifnum\pdfoutput<1 \else\pdftrue\fi\fi\fi
\ifpdf
\usepackage[hidelinks,bookmarks,bookmarksopen,bookmarksnumbered,pdfauthor={Daniel Berwick-Evans, Dmitri Pavlov},pdftitle={Smooth one-dimensional topological field theories are vector bundles with connection},backref=section,linktoc=all]{hyperref}
\usepackage[hiresbb]{graphics}
\else
\usepackage[hypertex,backref=section,linktoc=all]{hyperref}
\usepackage[hiresbb,dvips]{graphics}
\fi

\def\infinity{\texorpdfstring{$\infty$}{∞}}

\usepackage[capitalize]{cleveref}

\def\typeout#1{} 

\theoremstyle{plain}
\newtheorem{thm}{Theorem}[section]
\newtheorem{prop}[equation]{Proposition}

\newtheorem{cor}[equation]{Corollary}
\newtheorem{lem}[equation]{Lemma}

\theoremstyle{definition}
\newtheorem{defn}[equation]{Definition}

\theoremstyle{remark}
\newtheorem{rmk}[equation]{Remark}

\def\N{{\bf N}}

\def\R{{\bf R}}

\def\Vect{{\sf Vect}}
\def\Stacks{{\sf Stacks}}
\def\PreStacks{{\sf PreStacks}}
\def\Bord{\hbox{-{\sf Bord}}}
\def\CCat{\hbox{\sf$\sf C^\infty$-Cat}}

\def\sSet{{\sf sSet}}
\def\Set{{\sf Set}}
\def\CSS{{\hbox{\sf$\sf \infty$-Cat}}}

\def\Cart{{\sf Cart}}
\def\Cat{{\sf Cat}}
\def\SymCat{{\sf SymCat}}

\def\Mfld{{\sf Mfld}}

\def\op{{\sf op}}
\def\Hom{{\sf Hom}}
\def\Grpd{{\sf Grpd}}

\def\C{{\cal C}}
\def\D{{\cal D}}
\def\V{{\cal V}}

\def\path{{\cal P}}

\def\holim{\mathop{\rm holim}\nolimits}
\def\hocolim{\mathop{\rm hocolim}\nolimits}
\def\colim{\mathop{\rm colim}\nolimits}
\def\Aut{\mathop{\rm Aut}}
\def\End{\mathop{\rm End}}
\def\Map{\mathop{\rm Map}}

\def\Fun{\mathop{\rm Fun}}

\def\Ad{{\rm Ad}}
\def\Tran{{\rm Tran}}

\def\B{{\rm B}}
\def\O{{\rm O}}
\def\GL{{\rm GL}}

\def\pt{{\rm pt}}

\def\fr{{\rm or}}

\def\id{{\rm id}}
\def\equiv{{\rm equiv}}
\def\ci{C\/$^\infty$}
\def\Ci{{\rm C^\infty}}
\def\Citensor{{\rm C^\infty_\otimes}}

\def\To{\longrightarrow}

\def\[#1]{\langle#1\rangle}

\newcommand\tripleleftarrow{
\mathrel{\substack{\textstyle\leftarrow\\[-0.5ex]
                      \textstyle\leftarrow \\[-0.5ex]
                      \textstyle\leftarrow}}
}

\newcommand\quadleftarrow{
\mathrel{\substack{\textstyle\leftarrow\\[-0.5ex]
                      \textstyle\leftarrow \\[-0.5ex]
                      \textstyle\leftarrow \\[-0.5ex]
                      \textstyle\leftarrow}}
}

\def\rx#1#2{\rlap{\kern #1pt \raise#1pt \hbox{#2}}}
\def\dottednearrow{\rx{-8}. \rx{-6}. \rx{-4}. \rx{-2}. \rx0. \rx2. \rx4. \kern6pt \raise7.7pt \hbox{$\nearrow$}}

\catcode`@=11
\def\cal{\fam\tw@}
\let\over\@@over
\let\atop\@@atop
\let\above\@@above
\let\overwithdelims\@@overwithdelims
\let\atopwithdelims\@@atopwithdelims
\let\abovewithdelims\@@abovewithdelims
\newskip\xcentering \global\xcentering=0pt plus 1000pt minus 1000pt
\def\eqalign#1{\null\,\vcenter{\openup\jot\m@th
  \ialign{\strut\hfil$\displaystyle{##}$&$\displaystyle{{}##}$\hfil
      \crcr#1\crcr}}\,}
\def\eqalignno#1{\displ@y \tabskip\xcentering
  \halign to\displaywidth{\hfil$\@lign\displaystyle{##}$\tabskip\z@skip
    &$\@lign\displaystyle{{}##}$\hfil\tabskip\xcentering
    &\llap{$\@lign##$}\tabskip\z@skip\crcr
    #1\crcr}}
\def\cases#1{\left\{\,\vcenter{\normalbaselines\m@th
    \ialign{$##\hfil$&\quad##\hfil\crcr#1\crcr}}\right.}
\def\eqlabel#1{\refstepcounter{equation}\leqno(\theequation)\label{#1}}
\catcode`@=12

\def\cdiag#1{\null\,\vcenter{\normalbaselines
    \ialign{\hfil$##$\hfil&&\;\;\hfil$##$\hfil\crcr
      \mathstrut\crcr\noalign{\kern-\baselineskip}
      #1\crcr\mathstrut\crcr\noalign{\kern-\baselineskip}}}\,}
\def\cd{\def\normalbaselines{\baselineskip20pt \lineskip1pt \lineskiplimit0pt }\cdiag}

\def\mapdown#1{\Big\downarrow\rlap{$\vcenter{\hbox{$\scriptstyle#1$}}$}}

\def\matrix#1{\null\,\vcenter{\normalbaselines
        \ialign{\hfil$##$\hfil&&\enspace\hfil$##$\hfil\crcr
                \mathstrut\crcr\noalign{\kern-\baselineskip}
                #1\crcr\mathstrut\crcr\noalign{\kern-\baselineskip}}}\,}
\def\sqmatrix#1{\null\,\vcenter{\normalbaselines
        \ialign{\hfil$##$&\enspace\hfil$##$\hfil\thinspace&$##$\hfil\crcr
                \mathstrut\crcr\noalign{\kern-\baselineskip}
                #1\crcr\mathstrut\crcr\noalign{\kern-\baselineskip}}}\,}
\def\cdbl{\def\normalbaselines{\baselineskip20pt \lineskip3pt \lineskiplimit3pt }}
\def\cd{\cdbl\let\dskip\relax\matrix}
\def\sqcd{\cdbl\let\dskip\enskip\sqmatrix}

\newcount\arrowsize \arrowsize3
\def\Toarr#1{\mathop{\count0=#1 \loop\ifnum\count0>0 \smash-\mkern-7mu \advance\count0 -1 \repeat \mathord\rightarrow}\limits}
\def\xTo#1#2{\mathrel{\Toarr{#1}^{#2}}}
\def\Getsarr#1{\mathop{\mathord\leftarrow \count0=#1 \loop\ifnum\count0>0 \mkern-7mu\smash-\advance\count0 -1 \repeat}\limits}

\def\mapdown#1{\dskip\Big\downarrow\rlap{$\vcenter{\hbox{$\scriptstyle#1$}}$}\dskip}

\newcommand\triplerightarrow{
\mathrel{\substack{\textstyle\rightarrow\\[-0.5ex]
                      \textstyle\rightarrow \\[-0.5ex]
                      \textstyle\rightarrow}}
}

\newcommand\quadrightarrow{
\mathrel{\substack{\textstyle\rightarrow\\[-0.5ex]
                      \textstyle\rightarrow \\[-0.5ex]
                      \textstyle\rightarrow \\[-0.5ex]
                      \textstyle\rightarrow}}
}

\begin{document}

\title{Smooth one-dimensional topological field theories \\ are vector bundles with connection}

\author{Daniel Berwick-Evans\\
Department of Mathematics, University of Illinois Urbana--Champaign\\
\url{https://danbe.web.illinois.edu/}
\and
Dmitri Pavlov\\
Department of Mathematics and Statistics, Texas Tech University\\
\url{https://dmitripavlov.org/}}

\date{}

\maketitle

%$$

\begin{abstract}
We prove that smooth 1-dimensional topological field theories over a manifold are equivalent to vector bundles with connection.
The main novelty is our definition of the smooth 1-dimensional bordism category, which encodes \emph{cutting} laws rather than \emph{gluing} laws.
We make this idea precise through a smooth version of Rezk's complete Segal spaces.
With such a definition in hand, we analyze the category of field theories using a combination of descent,
a smooth version of the 1-dimensional cobordism hypothesis, and standard differential-geometric arguments.
\end{abstract}

\tableofcontents

\section{Introduction}

\numberwithin{equation}{section}

The goal of this paper is to give a definition of smooth 1-dimensional field theory that plays well with the differential geometry of manifolds.
A technical ingredient in our approach is the notion of \emph{smooth $\infty$-categories} (developed in \cref{sec:model}).
A smooth $\infty$-category is a smooth version of a complete Segal space. This framework can be used to encode cutting axioms for the value of a field theory on a cobordism rather than the usual gluing axioms.
This perspective on cobordisms is implicit in the work of Galatius--Madsen--Tillman--Weiss~\cite{GalatiusMadsenTillmannWeiss} and Lurie~\cite{Lurie:TFT} in the settings of topological and $\infty$-categories, respectively.
Translating into the smooth setting yields field theories that are determined by familiar differential-geometric objects.

{\def\thethm{A}
\begin{thm}
\label{thma}
The space of smooth 1-dimensional oriented topological field theories over a smooth manifold~$X$ is equivalent to the nerve of the groupoid of (finite-dimensional) vector
bundles with connection over~$X$ and connection-preserving vector bundle isomorphisms.
The equivalence is natural in~$X$.
\end{thm}
}

In our view, the above characterization of smooth 1-dimensional topological field theories is the only admissible one.
As such, the main contribution of this paper is a precise definition of smooth field theory for which \cref{thma} holds.
The definition readily generalizes both to higher dimensions and nontopological smooth field theories, as pursued by Grady and the second author~\cite[\S4]{GradyPavlov.Local}.
Through its connection to familiar objects, \cref{thma} gives a concrete idea of what these more complicated field theories seek to generalize.
Our methods---particularly the role of descent in \cref{thm:local}---are chosen with higher-dimensional generalizations in mind, see~\cite[Theorem 1.0.1]{GradyPavlov.Local}.

The intuition behind \cref{thma} goes back to Segal~\cite[Section~6]{Segal:EC}.
The 1-dimensional bordism category over~$X$ has objects compact 0-manifolds with a map to~$X$ and morphisms compact 1-manifolds with boundary with a map to~$X$.
A 1-dimensional topological field theory over~$X$ is a symmetric monoidal functor from the 1-dimensional bordism category over~$X$ to the category of vector spaces.
Hence, to each point in~$X$ a topological field theory assigns a finite-dimensional vector space and to each path the field theory assigns a linear map.
A vector bundle with connection produces this data using parallel transport.
Conversely, given the data of parallel transport one may assign values to arbitrary 1-dimensional bordisms in~$X$ by taking tensor products and duals. 
This gives a functor from the groupoid of vector bundles with connection to the groupoid of 1-dimensional topological field theories. 
To verify \cref{thma}, one must show that this functor is an equivalence. This can be thought of as the following smoothly-parametrized variant of the 1-dimensional cobordism hypothesis of Baez--Dolan~\cite{BaezDolan}. We note that in dimension~1, orientations are equivalent to framings.

{\def\thethm{B}
\begin{thm}
\label{thmb}
Let $\Vect^\otimes$ denote the symmetric monoidal smooth $\infty$-category of vector spaces and $\Vect$ the underlying smooth $\infty$-category without monoidal structure.
There is an equivalence between 1-dimensional oriented topological field theories over~$X$ valued in~$\Vect^\otimes$
and \ci-functors from the smooth path category of~$X$ to~$\Vect$.
\end{thm}
}

\cref{thma} follows from \cref{thmb} by identifying a functor from the path category to~$\Vect$ with a smooth vector bundle and connection (see \cref{sec:paths}).
Such a relationship between parallel transport and representations of path categories goes back to Kobayashi~\cite{Kobayashi},
who introduced the group of smooth based loops modulo thin homotopy
and established that smooth homomorphisms from this group to a Lie group~$G$
are in bijection with isomorphism classes of principal $G$-bundles with connection.
Similar results were proved by Freed~\cite{Freed}, and Schreiber--Waldorf~\cite{SchreiberWaldorf}.
An analogous result for gerbes with connections was established for two-dimensional thin homotopies by Bunke--Turner--Willerton~\cite{BunkeTurnerWillerton}.
A statement explicitly relating field theories to vector bundles with connection was loosely formulated by Segal in his early work on geometric models for elliptic cohomology.
A precise statement was given by Stolz and Teichner in their language of field theories fibered over manifolds~\cite[Theorem 1.8]{StolzTeichner.SUSY}, though a proof has yet to appear. 
Below we draw inspiration from all of these authors;
the new ingredient is in our treatment of the smooth bordism category.

\numberwithin{equation}{subsection}

\subsection{What makes a smooth cobordism category difficult to define?}

Composition of cobordisms is more subtle than one might hope:
given two $d$-manifolds and a $(d-1)$-manifold along which one wishes to glue, one only obtains a glued manifold \emph{up to diffeomorphism}.
The usual solution is to define morphisms in a cobordism category as smooth $d$-manifolds up to diffeomorphism, thereby obliterating the problem.
However, if we wish to include extra structures on bordisms, two problems arise:
(1)~gluings may fail to exist and (2)~gluing isomorphism classes of geometric structures is typically ill-defined. 
For example, 
isomorphism classes of metrics cannot be glued, and gluing smooth maps to a target manifold~$X$ along a codimension~1 submanifold may not result in a smooth map to~$X$.
Although the focus in this paper is on the topological bordism category over~$X$, a guiding principle is to make definitions for which various flavors of generalization pose no serious technical difficulties.

A reason for pursuing such generalizations comes from Stolz and Teichner's work on a geometric model for elliptic cohomology~\cite{StolzTeichner.Elliptic, StolzTeichner.SUSY}, following Segal~\cite{Segal:EC}.
Their program seeks to generalize the relationship between 1-dimensional field theories, vector bundles, and K-theory to provide a model for elliptic cohomology with cocycles coming from (supersymmetric) 2-dimensional Euclidean field theories. 
They address the problem of composition by equipping bordisms with germs of collars, so bordisms are composable when collars match.
This allows one to incorporate geometric structures on paths by simply endowing the collars with geometric structures.
However, it introduces a new technical issue when relating 1-dimensional field theories to vector bundles: isomorphism classes of objects in their 1-dimensional bordism category are points of~$X$ \emph{together} with the germ of a collar of a path.
Such a large space of objects can be rather unwieldy in computations (particularly in the presence of geometric structures on the bordisms). The original motivation for this paper was to find a way around the technical difficulties brought on by the introduction of collars, with an eye towards studying supersymmetric Euclidean field theories. 
We note that since the writing of this paper, an analog of \cref{thma} was proved by Ludewig and Stoffel~\cite{LudewigStoffel} in a framework that incorporates collared bordisms (we comment further on this work in~\cref{sec:subsequent}).

One approach that avoids collars follows Kobayashi~\cite{Kobayashi}, Caetano--Picken~\cite{CaetanoPicken}, and Schreiber--Waldorf \cite{SchreiberWaldorf} who study paths in a manifold modulo \emph{thin homotopy};
these are smooth paths modulo smooth homotopies whose rank is at most~1.
Each equivalence class has a representative given by a path with sitting instant, meaning a path in~$X$ for which
some neighborhood of the start and end point is mapped constantly to~$X$.
These sitting instants allow concatenation of smooth paths in a straightforward way, which simplifies many technical challenges (compare \cref{lem:sitting}), and leads to a version of a 1-dimensional bordism category over~$X$.
However, endowing an equivalence class of a path with a geometric structure, e.g., a metric, is hopeless.
If one does not pass to equivalence classes but instead works with honest paths with sitting instants, the resulting path category fails to restrict along open covers of~$X$: restricted paths may not have sitting instants.
This destroys a type of locality that we find both philosophically desirable and computationally essential (compare \cref{thm:local}), related to Mayer--Vietoris sequences in Stolz and Teichner's program, see~\cite[Conjecture 1.17]{StolzTeichner.SUSY}.
In summary, techniques involving bordisms with sitting instants seem appropriate only for a certain class of topological field theories.

Lurie~\cite{Lurie:TFT} and Galatius--Madsen--Tillman--Weiss~\cite{GalatiusMadsenTillmannWeiss} take a different road, considering a topological category (or Segal space) of bordisms, wherein composition need only be defined up to homotopy.
This allows one to effectively add or discard the collars with impunity since this data is contractible.
Geometric structures can be incorporated as stable tangential structures, i.e., maps from bordisms to certain classifying spaces.
This framework also leads to field theories that are relatively easy to work with:
one can obtain a precise relationship between maps from $X$ to~$\B\O(k)$ and 1-dimensional topological field theories over~$X$.
However, the price one pays is that such bundles are not smooth, but merely topological.
By this we mean a particular space of field theories is homotopy equivalent to the space of maps from~$X$ (viewed as a topological space, not a manifold)
to the classifying space of vector bundles; in particular, from this vantage the data of a connection is contractible.
Our search for a smooth bordism category is tantamount to asking for a differential refinement of this data.
In the case of line bundles, such a refinement is the jumping-off point for the subject of (ordinary) differential cohomology,
and one can view our undertaking as a close cousin.

\subsection{Why model categories?}
Our approach combines aspects of Stolz and Teichner's definition of bordism categories internal to smooth stacks \cite{StolzTeichner.Elliptic, StolzTeichner.SUSY}
and the Segal space version (in the world of model categories) studied by Lurie.
What we obtain is not a bordism \emph{category} strictly speaking, but rather bordisms in~$X$ comprise a collection of objects and morphisms with
a partially defined composition; this is the categorical translation of the geometric idea to encode cutting laws rather than gluing laws.
To make sense of functors out of this bordism ``category'' we provide ourself with an ambient category of smooth
categories with partially defined composition. 
This is directly analogous to Rezk's category of complete Segal spaces as a model for $\infty$-categories.
However, making the framework precise requires a foray into the world of model categories.
This sort of machinery rarely turns up in standard differential geometry, so might seem a little misplaced at first glance; however,  in the setting of field theories some basic features of this language seem unavoidable.
For example, the bordism category over~$X$ ought to be equivalent to the bordism category over an open cover $\{U_i\to X\}$
with appropriate compatibility conditions on intersections; asking for
these categories to be \emph{isomorphic} is too strong since, for example, they have different sets of objects.
Hence, the appropriate categorical setting for describing bordisms over~$X$ must have some native notion of (weak) equivalence of bordism categories,
and the language of model categories was built precisely to facilitate computations in such situations.

Smooth $\infty$-categories are fibrant objects for a model category structure we place on simplicial objects in smooth stacks.
This model structure is chosen to satisfy three properties:
(1)~nerves of categories fibered over manifolds determine fibrant objects, e.g., the category of (smooth) vector spaces is fibrant;
(2)~all objects are cofibrant, and in particular bordisms over~$X$ define a cofibrant object; and
(3)~there is a weak equivalence between 1-dimensional bordisms over~$X$ and 1-dimensional bordisms over an open cover~$\{U_i\}$ of~$X$ with compatibility on overlaps.
These three properties are the only features of the model structure we actually use.
It is important that bordisms over~$X$ do not give a fibrant object: this is precisely the failure of composition to be defined in general.
However, to compute the (derived mapping) space of smooth functors from bordisms to vector spaces, fibrant replacement of the source is unnecessary.

We refer the reader to Hirschhorn \cite{Hirschhorn}, Barwick \cite{Barwick:LR}, and Lurie \cite[Appendix~A]{Lurie:HTT} for background on model categories relevant to this paper.

\subsection{Subsequent work}\label{sec:subsequent}

Since the first version of this paper appeared in 2015, its basic ideas have been used and expanded by several authors. Some of these developments were alluded to above; we discuss this further presently. 

Ludewig and Stoffel~\cite{LudewigStoffel} constructed a bordism category using model-categorical techniques, incorporating ideas similar to those in this paper.
One important difference is that their bordisms are equipped with the germ of a collar, following the ideas of Stolz and Teichner~\cite{StolzTeichner.SUSY}.
Ludewig--Stoffel go on to prove a version of \cref{thma}, see~\cite[Theorems 1.1 and 1.2]{LudewigStoffel}.
This shows that incorporating collars in the context of 1-dimensional topological field theories has no effect on the underlying geometric objects of study. 
Another important result of Ludewig--Stoffel shows that a version of \cref{thma} holds where the target category consists of (not necessarily locally free) sheaves of vector spaces~\cite[Theorem 5.2]{LudewigStoffel}.
These more general sheaves are an important piece of Stolz and Teichner's formalism, see~\cite[Remark 3.16]{StolzTeichner.SUSY}.
As in the formalism below, a crucial tool in Ludewig and Stoffel's work is the descent property for field theories.
We also mention that the main results of Benini--Perin--Schenkel~\cite{BeniniPerinSchenkel} verify and utilize descent for a distinct (though related) category of 1-dimensional \emph{algebraic} quantum field theories.

In~\cite{GradyPavlov.Local}, Grady and the second author generalize the 1-dimensional topological bordism category
defined in this paper to arbitrary $d$-dimensional bordism $d$-categories.
Their definitions also allow one to incorporate a wide class of geometric structures on $d$-dimensional bordisms,
expanding the ideas of Stolz and Teichner~\cite{StolzTeichner.SUSY} to the fully-extended context.
The categorical foundations of Grady--Pavlov involve a smooth refinement of $d$-fold iterated Segal spaces,
generalizing the definition of smooth $\infty$-category developed below.
The main result of Grady--Pavlov~\cite{GradyPavlov.Local} is that fully extended $d$-dimensional geometric field theories satisfy descent;
their proof generalizes the ideas of \cref{sec:model}.
This descent statement leads to good behavior for groupoids of field theories over a manifold.
In particular, concordance classes of $d$-dimensional geometric field theories are representable~\cite[\S7.2]{GradyPavlov.Local}.
This fits nicely with Stolz and Teichner's point of view: the main conjectures from~\cite[\S1]{StolzTeichner.SUSY}
together with Brown representability require that concordance classes of $1|1$- and $2|1$-dimensional Euclidean field theories be representable.
Descent gives a conceptually satisfying mechanism for representability of a wide class of geometric field theories. 

Another important development is a proof (by Grady and the second author \cite{GradyPavlov.GCH}) of a geometric version of the cobordism hypothesis.
The model categorical framework is crucial, as freely adjoining duals to a smooth $(\infty,d)$-category can be realized by Bousfield localization.
The formalism and techniques of \cref{sec:model} continue to play a central role: generalizations of cutting-and-gluing constructions decompose bordisms into elementary handles that in turn generate the bordism categories of interest. 

\subsection{Notation and terminology}

\begin{defn}
\label{cartesian.site}
Let $\Cart$ denote the \emph{cartesian site} whose objects are~$\R^n$ for~$n\in\N$,
morphisms are all smooth maps, and coverings are the usual open coverings, $\coprod\R^n\to\R^n$.
\end{defn}

We will use the notation~$[k]$ to denote the finite set~$\{0,1,\ldots,k\}$ as an object of the category $\Delta$ of simplices. The word `space' will often be used to refer to simplicial sets (i.e., objects in the category $\sSet=\Fun(\Delta^\op,\Set)$), e.g., a simplicial space is a bisimplicial set.

We will sometimes refer to objects in $\CCat$ and $\CCat^\otimes$ as categories even when they are not, e.g., we will often refer to $1\Bord^\fr(X)$
as the 1-dimensional oriented bordism category over~$X$.

\subsection{Structure of the paper}

We have attempted to keep the model-categorical discussion separate from the geometric one,
relegating the former to \cref{sec:model} and the latter to Sections~\ref{sec:1bord}--\ref{sec:1TFT}.

\subsection{Acknowledgments}

We thank Peter Teichner for many enlightening discussions about field theories, 
Ralph Cohen for his encouragement, 
Konrad Waldorf for helpful correspondence about representations of path groupoids,
and Gabriel Drummond-Cole for feedback on an earlier draft.
Pavlov was partially supported by the SFB 878 grant.

\section{From \infinity-categories to smooth \infinity-categories}
\label{sec:model}

Complete Segal spaces provide a model for $\infty$-categories.
A complete Segal space 
is a simplicial space (i.e., a functor $\Delta^\op\to \sSet$) satisfying a Segal and completeness condition reviewed below. 
As described by Rezk~\cite{Rezk}, the value of a Segal space on $[k]\in \Delta$ can be thought of as the classifying space for chains of morphisms of length~$k$ in an ordinary category.
This has an obvious smooth enhancement via simplicial objects in smooth stacks~$\Delta^\op\to \Stacks$ wherein the value on $[k]$ is a classifying stack for chains of composable morphisms of length~$k$.
This is the approach we take leading to the definition of a smooth $\infty$-category.
By adjunction, one can also view these smooth $\infty$-categories as sheaves of complete Segal spaces on the smooth site~$\Cart$. 

For our applications to field theories, we will assemble 1-dimensional bordisms in~$M$ into a functor $\Delta^\op\to {\rm Stacks}$
whose value on $[k]\in \Delta$ is the classifying stack of (smoothly) composable chains of bordisms in $M$ of length~$k$.
As we shall see, this functor does \emph{not} satisfy the Segal condition.
The problem is geometric and unavoidable: arbitrary chains of bordisms in~$M$ compose to a piecewise smooth bordism that need not be smooth. 
To work with such an object, we require a larger category that includes simplicial objects in stacks that do not satisfy the Segal condition. 

A systematic method for dealing with this type of issue is to construct a model category whose fibrant objects satisfy a Segal and completeness condition.
This allows one to work with nonfibrant objects precisely when their failure to be fibrant is not homotopically problematic. 
For example, mapping out of a nonfibrant object usually presents no issues, whereas mapping into a nonfibrant object can be problematic. 
For complete Segal spaces, such a model structure was constructed by Rezk as a localization of the Reedy model structure on simplicial spaces.
It is only a mild elaboration to extend these ideas to simplicial objects in stacks, i.e., smooth $\infty$-categories.

\subsection{Complete Segal spaces as fibrant objects in a model category}
We overview the small part of Rezk's theory of complete Segal spaces that we require; see Rezk~\cite{Rezk} for a more thorough treatment. 

\begin{defn}
\label{segal.space}
A \emph{Segal space} is a functor $C\colon\Delta^\op\to\sSet$ that satisfies the Segal condition, meaning the \emph{Segal map}
$$C(k)\to C(1)\times_{C(0)}^h\cdots\times_{C(0)}^hC(1),\quad k\ge 1\eqlabel{eq:Segmap}$$
into the homotopy fibered product is a weak equivalence.
\end{defn}

When the spaces are Reedy fibrant, the homotopy fibered product can be computed as the ordinary fibered product.

To obtain a category of Segal spaces into which the category of categories can be naturally embedded, it is important to also enforce a \emph{completeness condition} that we recall presently.
For any Segal space~$C$, let $h C$ denote the underlying homotopy category.
This is an ordinary category defined by Rezk~\cite[\S5]{Rezk},
whose objects are 0-simplices of~$C(0)$
and whose morphisms from~$x_0$ to~$x_1$ are the connected components of the homotopy fiber of~$C(1)$ over~$(x_0,x_1)$
for the projection $$d_0\times d_1\colon C(1)\to C(0)\times C(0).$$
For a Segal space~$C$, let $C_\equiv\subset C(1)$ be the subspace consisting of connected components of the above fiber that correspond to isomorphisms in the homotopy category.
Rezk shows that the degeneracy map $s_0\colon C(0)\to C(1)$ factors through the subspace~$C_\equiv$. 

\begin{defn}
\label{complete.segal.space}
A Segal space (\cref{segal.space}) is \emph{complete} if the map 
$$s_0\colon C(0)\to C_\equiv\eqlabel{eq:completeness}$$
is a weak equivalence.
\end{defn}

Following Rezk~\cite[\S7]{Rezk}, we can define a model category $\CSS$ whose underlying category is simplicial spaces and whose fibrant objects are complete Segal spaces.
We achieve this by localizing a given standard model structure,
namely, the Reedy model structure on simplicial spaces,
which coincides with the injective model structure, see, for example, Bergner--Rezk \cite[Proposition 3.15 (arXiv); 3.10 (journal); Proposition 4.2 (arXiv); 4.1 (journal)]{BergnerRezk}.
For a brief review of Reedy model structures, see \cref{Reedy}.
The first set of morphisms in this localization come from the morphisms of simplicial presheaves on~$\Delta$ for each $k\in \N$ given by
$$\varphi_k\colon [1]\sqcup_{[0]}^h[1]\sqcup_{[0]}^h\cdots \sqcup_{[0]}^h [1]\to [k]\eqlabel{eq:Segalsp}$$
where the source is a $k$-fold iterated homotopy pushout in the category of simplicial presheaves on~$\Delta$ (alias: simplicial spaces),
and we have identified $[m]\in \Delta$ with its associated (representable) presheaf.
The maps $[0]\to [1]$ defining the homotopy pushout $[1]\sqcup_{[0]}^h [1]$ are $0\mapsto 0$ and $0\mapsto 1$, and the higher pushouts iterate this basic version.
The homotopy colimit in~\eqref{eq:Segalsp} can be computed as the ordinary colimit because the underlying diagram is cofibrant.
If we map the source and target of~$\varphi_k$ into a simplicial space~$C$ and consider the map that $\varphi_k$ induces on these mapping spaces, we obtain the Segal map~\eqref{eq:Segmap}. 
(The resulting map $\varphi_k$ has a cofibrant domain and codomain, which means that derived mapping spaces out of them can be computed as ordinary mapping spaces.)

Let $E$ denote the simplicial space associated to the simplicial set given by the nerve of the groupoid with two objects $x$ and~$y$,
and two nonidentity morphisms $x\to y$ and $y\to x$.
Consider the canonical map
$$x\colon E\to [0],\eqlabel{eq:completesp}$$ 
where again we have identified $[0]\in \Delta$ with its associated representable simplicial presheaf.
Rezk~\cite[Theorem~6.2]{Rezk} shows that $\Map(E,C)\simeq C_\equiv$ and the map 
$$C(0)\to \Map(E,C)\simeq C_\equiv$$
induced by~\eqref{eq:completesp} is a weak equivalence if and only if $C$ is complete.
(In fact, this map is weakly equivalent to the map \eqref{eq:completeness} defined above.)

\begin{defn}
Endow the category $\Fun(\Delta^\op,\sSet)$ of simplicial spaces with the Reedy model structure, which coincides with the injective structure.
Define the \emph{model category of complete Segal spaces}, denoted by $\CSS$,
as the left Bousfield localization of this model structure along the maps $\varphi_k$ and $x$ from~\eqref{eq:Segalsp} and~\eqref{eq:completesp}.
\end{defn}

\begin{rmk}
\label{fibrant.css}
It is immediate from the properties of the localized model structure that fibrant objects in $\CSS$ are simplicial spaces that are fibrant in the Reedy model structure and satisfy the Segal and completeness conditions; see Rezk~\cite[Theorem~7.2]{Rezk}, or compare the proof of \cref{fibrant.smooth.cat} below. 
In particular, fibrant objects in $\CSS$ coincide with Reedy fibrant complete Segal spaces in the sense of \cref{complete.segal.space}.
\end{rmk}

\subsection{Smooth \infinity-categories}

As mentioned at the beginning of the section, we take smooth $\infty$-categories to be a stack-valued version of Segal spaces.
A \emph{(smooth) stack} is a functor~$F\colon \Cart^\op\to \sSet$ satisfying descent for good open covers~$\{U_i\}$ of objects $S\in \Cart$, meaning the canonical map
$$F(S)\longrightarrow \holim\left(\prod F(U_i)\rightrightarrows \prod F(U_{ij}) \triplerightarrow \prod F(U_{ijk}) \quadrightarrow\cdots\right)\eqlabel{homotopy.descent}$$
is a weak equivalence.
It will be useful later to observe this comes from mapping the source and target of 
$$S\stackrel{p}{\longleftarrow} \hocolim\left(\coprod_i U_i\leftleftarrows \coprod_{i,j} U_i\cap U_j \tripleleftarrow \coprod_{i,j,k} U_i\cap U_j\cap U_k \quadleftarrow\cdots\right)\eqlabel{eq:Cech}$$
to~$F$ and considering the map induced by~$p$ between the resulting mapping spaces of simplicial presheaves on~$\Cart$.
In this description, the above is regarded as a morphism of simplicial presheaves and (in particular) we have identified~$S$
and the $U_i$ with their representable presheaves and taken the homotopy colimit in simplicial presheaves.
The latter homotopy colimit can be computed as the \v Cech nerve of~$U$ in simplicial presheaves (where $U_i\cap U_j\cap U_k$ is placed in simplicial degree~2,
and analogously for other intersections)
due to the Reedy cofibrancy of the underlying simplicial diagram.
The resulting morphism has a projectively cofibrant domain and codomain, which allows us to compute derived mapping spaces out of them as ordinary mapping spaces
if the target is projectively fibrant, i.e., an objectwise Kan complex.
We emphasize that in~\eqref{eq:Cech} the coproduct $\coprod$ is taken in the category of presheaves (which is different from the category of sheaves).
The original reference for the model structure on stacks is Jardine~\cite{Jardine}.
A description in terms of left Bousfield localizations can be found in Dugger--Hollander--Isaksen~\cite{DuggerHollanderIsaksen}.

\begin{defn}
\label{model.stack}
The model category $\PreStacks$ is the projective model structure on simplicial presheaves on the cartesian site $\Cart$ (\cref{cartesian.site}), i.e., $\Fun(\Cart^\op,\sSet)$.
The model category $\Stacks$ is the left Bousfield localization of $\PreStacks$ along the morphisms~\eqref{eq:Cech} for all good open covers $\{U_i\}_{i\in I}$ of any $S\in\Cart$.
\end{defn}

\begin{rmk}
\label{fibrant.stack}
Fibrant objects in $\PreStacks$ (\cref{model.stack}) are precisely presheaves valued in Kan complexes,
whereas fibrant objects in $\Stacks$ are precisely those fibrant objects in $\PreStacks$ that satisfy the homotopy descent condition~\eqref{homotopy.descent}.
\end{rmk}

We now consider the Reedy model structure on the category of simplicial prestacks, i.e., functors $\Delta^\op\to \PreStacks$.
The existence and basic properties of this model structure follows from Hirschhorn~\cite[Theorem~15.3.4]{Hirschhorn}, as we review briefly.
Weak equivalences are objectwise, meaning a map $F\to G$ of simplicial prestacks is an equivalence if we get an equivalence of prestacks for each fixed $[n]\in \Delta$.
Fibrations and cofibrations are described in terms of relative matching and latching maps; see \cref{Reedy}.
By adjunction, the Reedy model structure also gives a model structure on the presheaf category
$$\Delta^\op\times \Cart^\op\to \sSet.$$
In this description, fibrant objects $C\colon \Delta^\op\times \Cart^\op\to \sSet$ are precisely those presheaves that define Reedy fibrant simplicial spaces~$C(-,S)\colon \Delta^\op\to \sSet$ for each~$S\in \Cart$. 

We wish to localize this Reedy model structure on simplicial prestacks along the morphisms~\eqref{eq:Segalsp}, \eqref{eq:completesp} and~\eqref{eq:Cech}, but to do so we need to promote these to morphisms in presheaves on $\Delta\times \Cart$.
We achieve this in the following way. 
For each representable presheaf associated to an object $[n]\in \Delta$ or $S\in \Cart$, we take the morphisms of presheaves
$$\id_{[n]}\times p,\quad \varphi_k\times \id_S, \quad x\times \id_S,\eqlabel{eq:infinityloc}$$ 
where $\varphi_k$, $x$ and $p$ are as in~\eqref{eq:Segalsp}, \eqref{eq:completesp} and~\eqref{eq:Cech} and $\id_{[n]}$ or $\id_S$ denotes the identity morphism on the corresponding representable presheaf. 
Letting $[n]\in \Delta$ and $S\in \Cart$ range over all possible objects, we obtain our localizing morphisms. 

\begin{defn}
\label{model.smooth.cat}
Define the \emph{model category of smooth $\infty$-categories}, denoted by $\CCat$,
as the left Bousfield localization of the Reedy model structure on $$\Fun(\Delta^\op,\PreStacks)$$ along the set of morphisms~\eqref{eq:infinityloc}. 
A \emph{smooth functor} is a morphism in $\CCat$.
\end{defn}

Existence and basic properties of such a localization is shown by Barwick~\cite[Theorem~2.11 (arXiv); 4.7 (journal)]{Barwick:LR},
since the injective or projective model structure on the category of simplicial presheaves is left proper and combinatorial,
see, for example, Lurie \cite[§A.2.7, Proposition~A.2.8.2, Remark~A.2.8.4]{Lurie:HTT}.
We summarize what we require as follows.

\begin{lem}
\label{fibrant.smooth.cat}
Fibrant objects $C\in \CCat$ (\cref{model.smooth.cat}) are simplicial presheaves on $\Delta\times \Cart$ such that
\begin{enumerate}
\item for any~$S\in \Cart$, the restriction $C(-,S)\colon \Delta^\op\to \sSet$ is a fibrant complete Segal space (\cref{fibrant.css});
\item for any $[n]\in \Delta$, the restriction $C([n],-)\colon \Cart^\op\to \sSet$ is a fibrant smooth stack (\cref{fibrant.stack}).
\end{enumerate}
\end{lem}

\begin{proof}
By Barwick~\cite[Theorem~2.11 (arXiv); 4.7 (journal)]{Barwick:LR}, an object $C$ is fibrant in the local model structure if it is fibrant in the Reedy model structure (i.e., before localization) and has the additional property that 
for all maps $f\colon A\to B$ in~\eqref{eq:infinityloc}, the induced map of derived mapping spaces
$$\Map(B,C)\to\Map(A,C)\eqlabel{eq:local}$$
is a weak equivalence.
As observed above, $C$ being fibrant before localization reduces to $C(-,S)$ being a Reedy fibrant simplicial space for any~$S\in \Cart$.
We observe that the maps in~\eqref{eq:infinityloc} have cofibrant source and fibrant target, 
so the derived mapping spaces in~\eqref{eq:local} can be computed as the usual mapping spaces.
For the maps $\phi_k\times\id_S$ and $x\times\id_S$, we see that $C(-,S)$ must be a complete Segal space.
For the maps $\id_{[n]}\times p$ we see that $C([n],-)$ must be a stack.
\end{proof}

\subsection{Symmetric monoidal smooth \infinity-categories}

For our intended application to field theories, we also require a version of symmetric monoidal smooth $\infty$-categories.
Following ideas of Segal~\cite{Segal:CCT}, we implement this via the category~$\Gamma$, the opposite category of finite pointed sets.

Just like for~$\Delta$, presheaves on $\Gamma$ are first equipped with the Reedy model structure, which is then localized with respect to appropriate maps.
The category~$\Gamma$ has nontrivial automorphisms, so the usual notion of a Reedy model structure
must be generalized to accommodate this new setting,
resulting in the {\it strict model structure\/} of Bousfield--Friedlander \cite[§3]{BousfieldFriedlander}.
This approach was generalized by Berger--Moerdijk \cite{BergerMoerdijk}, resulting in the notion of a {\it generalized Reedy category\/}
and the associated model structure.

As explained in Segal~\cite{Segal:CCT},
given a $\Gamma$-object~$X$, we can think of $X_{\[n]}=X_n$ as the space of (formal) $n$-tuples of elements of some commutative monoid.
Here “formal” means that points of $X_n$ are not actual $n$-tuples, but rather have certain structure that makes them formally behave like ones.
Specifically, given a map of finite pointed sets $f\colon\[m]\to\[n]$, the associated map $X_f\colon X_{\[m]}\to X_{\[n]}$
should be thought of as multiplying the elements indexed by $f^{-1}\{j\}$ for each $j\in\[n]$ (the product of an empty family is the identity element),
and throwing away elements indexed by $f^{-1}\{*\}$.
Given a $\Gamma$-object~$X$, its $n$th latching map $L_n X\to X_n$ can be thought of as the subobject of $X_n$ comprising those (formal) $n$-tuples
where at least one element is the identity.
Given a $\Gamma$-object~$X$, its $n$th matching map $X_n\to M_n X$ can be thought of as sending a formal $n$-tuples in~$X_n$
to the compatible family of formal tuples given by multiplying two or more elements, or throwing away one or more elements.
A $\Gamma$-object~$X$ is cofibrant in the strict model structure if for all $n\in\Gamma$,
the latching map $L_n X\to X_n$ is a projective cofibration of objects equipped with an action of~$\Sigma_n$.
A $\Gamma$-object~$X$ is fibrant in the strict model structure if for all $n\in\Gamma$,
the matching map $X_n\to M_n X$ is a fibration.

To define the Reedy model structure on presheaves on $\Delta\times\Gamma$,
it suffices to observe that generalized Reedy categories are closed under finite products by Berger--Moerdijk \cite[§1]{BergerMoerdijk}
and both $\Delta^\op$ and $\Gamma^\op$ have generalized Reedy category structures
by Berger--Moerdijk \cite[Examples 1.8(a,b) (arXiv), 1.9(a,b) (journal)]{BergerMoerdijk}.
Thus, we can consider the (generalized) Reedy model structure on the category of functors
$$\Delta^\op\times \Gamma^\op\to \PreStacks.$$

Next, we turn attention to the morphisms used to define a local model structure.
As above, we define these morphisms in their adjoint form in the category of functors
$$\Delta^\op\times \Gamma^\op\times \Cart^\op \to \sSet.$$
Denote objects of $\Gamma$ by $\[m]=\{*,1,\ldots,m\}$, where $*$ is the basepoint.
Then consider the set of morphisms
$$u\colon \emptyset\to\[0], \qquad t^{m,n}\colon \[m]\sqcup_{\[0]} \[n]\to \[m+n],$$
where the maps $t^{m,n}$ are induced by a pair of maps of finite pointed sets (which yield morphisms in~$\Gamma$ in the opposite direction) $\[m]\leftarrow \[m+n]$ and $\[n]\leftarrow \[m+n]$. 
The first of these maps is the identity on the subset $\{1,\ldots,m\}$ and sends $m+1,\ldots,m+n$ to $*$.
The second of these maps uses the obvious bijection from $\{m+1,\dots, m+n\}$ to $\{1,\dots,n\}$ and sends the remainder to~$*$.

We then consider a set of localizing morphisms similar to~\eqref{eq:infinityloc}, only now we have to fix objects in a pair of the categories $\Delta$, $\Gamma$ and $\Cart$.
Explicitly, we take morphisms 
$$\id_{[n]}\times\id_{\langle m\rangle}\times p,\quad \varphi_k\times\id_{\langle m\rangle}\times \id_S,\quad
x\times\id_{\langle m\rangle}\times \id_S,\quad \id_{[n]}\times t^{m,n}\times \id_S,\quad \id_{[n]}\times u\times \id_S,\eqlabel{eq:symmonloc}$$
where $[n]\in \Delta$, $\langle m\rangle\in \Gamma$, and $S\in \Cart$ vary over all possible objects.

\begin{defn}
\label{model.smooth.mon.cat}
Define the \emph{model category of symmetric monoidal smooth $\infty$-categories}, denoted by $\CCat^\otimes$,
as the left Bousfield localization of the Reedy
model structure on the category of functors $\Delta^\op\times \Gamma^\op\to \PreStacks$
with respect to the morphisms~\eqref{eq:symmonloc}. 
\end{defn}

Again, existence of such a localization is shown by Barwick~\cite[Theorem~2.11 (arXiv); 4.7 (journal)]{Barwick:LR},
which also proves that fibrant objects in the localized structure
are precisely fibrant and local objects in the original model structure.

\begin{lem}
\label{fibrant.smooth.mon.cat}
Fibrant objects $C\in \CCat^\otimes$ (\cref{model.smooth.mon.cat}) are Reedy fibrant simplicial presheaves $C\colon \Delta^\op \times \Gamma^\op\times \Cart^\op\to \sSet$
(meaning the adjoint map $\Delta^\op\times\Gamma^\op \to \PreStacks$ is Reedy fibrant, where $\PreStacks$ is equipped with the projective model structure of \cref{model.stack})
such that
\begin{enumerate}
\item for any~$S\in \Cart$ and $\[m]\in \Gamma$, the restriction $C(-,\[m],S)\colon \Delta^\op\to \sSet$ is a complete Segal space (\cref{fibrant.css}); 
\item for any fixed $[n]\in \Delta$ and $\[m]\in \Gamma$, the restriction $C([n],\[m],-)\colon \Cart^\op\to \sSet$ is a smooth stack (\cref{fibrant.stack});
\item for any fixed $S\in\Cart$ and $[n]\in\Delta$, the restriction $C([n],-,S)\colon\Gamma^\op\to\sSet$ is a special $\Gamma$-space (Segal \cite[Definition~1.2]{Segal:CCT}; Bousfield--Friedlander \cite[§4]{BousfieldFriedlander}).
\end{enumerate}
\end{lem}

\begin{proof}
Fibrant objects in the local model structure are fibrant object in the Reedy model structure (selected by the first condition)
that additionally satisfy the locality condition~\eqref{eq:local}
with respect to the maps~\eqref{eq:symmonloc}.
For the maps $\phi_k\times\id_{\[m]}\times\id_S$ and $x\times\id_{\[m]}\times\id_S$ we see using the adjunction property that $C(-,\[m],S)$ must be a complete Segal space.
For the maps $\id_{[n]}\times\id_{\[m]}\times p$ we see that $C([n],[m],-)$ must be a stack.
For the maps $\id_{[n]}\times t^{m,n}\times \id_S$ and $\id_{[n]}\times u\times \id_S$ we see that $C([n],-,S)$ must be a special $\Gamma$-space.
\end{proof}

\begin{defn}
\label{forget.mon}
We have a functor
$$\CCat^\otimes \to \CCat$$ 
that restricts a simplicial presheaf on $\Delta\times \Gamma\times \Cart$ to $\Delta\times \langle 1\rangle\times \Cart\cong \Delta\times \Cart$.
We call this the \emph{forgetful functor} from symmetric monoidal smooth $\infty$-categories to smooth $\infty$-categories.
\end{defn}

By virtue of \cref{fibrant.smooth.cat} and \cref{fibrant.smooth.mon.cat},
the forgetful functor preserves fibrant objects and weak equivalences between them, and in fact is a right Quillen functor, which can be seen as follows.
Hirschhorn~\cite[Theorem~15.5.2]{Hirschhorn} shows that there is no difference between the Reedy model structures
on $$\Fun(\Delta^\op,\Fun(\Gamma^\op,\PreStacks)), \qquad \Fun(\Gamma^\op,\Fun(\Delta^\op,\PreStacks)), \qquad\Fun(\Delta^\op\times\Gamma^\op,\PreStacks),$$
using that $\Delta^\op\times\Gamma^\op$ is a (generalized) Reedy category,
with Hirschhorn's proof still working for generalized Reedy categories.
In particular, the forgetful functor can be presented as evaluation at $\[1]\in\Gamma$:
$$\Fun(\Gamma^\op,\Fun(\Delta^\op,\PreStacks))\to \Fun(\Delta^\op,\PreStacks),$$
hence it is a right Quillen functor because Reedy (acyclic) fibrations are projective (acyclic) fibrations.

\subsection{Smooth \infinity-categories from sheaves of categories}\label{sec:cifibered}
\label{examples.machine}

The main example of a Segal space comes from a simplicial space-valued nerve of an ordinary category; see Rezk~\cite[\S3.3]{Rezk}.
For small categories $\C$ and~$\D$, let ${\rm Iso}(\C^{\D})$ denote the category whose objects are functors $\D\to \C$ and whose morphisms are natural isomorphisms of functors.
Then define
$$\N^\infty\colon\Cat\to \Fun(\Delta^\op,\sSet),\qquad  \C\mapsto([l]\mapsto {\rm N}({\rm Iso}(\C^{[l]}))).\eqlabel{eq:cssnerve}$$
Rezk~\cite[Proposition~6.1]{Rezk} proves that $\N^\infty(C)$ is a complete Segal space; we sketch the idea.
The set of $(k,l)$-bisimplices of $\N^\infty(\C)$ is the set of diagrams in~$\C$
$${\rm Nerve}_k({\rm Iso}(\C^{[l]}))=\left\{\cd{c_{00}&\xTo2{f_{10}}&c_{10}&\xTo2{f_{20}} & \cdots & \xTo2{f_{l0}} & c_{l0}\cr
\mapdown{}&&\mapdown{}& & & & \mapdown{}  \cr
c_{01}&\smash{\xTo2{f_{11}}}& c_{11}& \xTo2{f_{21}} & \cdots & \xTo2{f_{l1}} & c_{l1}\cr
\mapdown{}&&\mapdown{}& & & & \mapdown{}  \cr
\vdots{}&&\vdots{}& & & & \vdots{}  \cr
\mapdown{}&&\mapdown{}& & & & \mapdown{}  \cr
c_{0k}&\smash{\xTo2{f_{1k}}}& c_{1k}& \xTo2{f_{2k}} & \cdots & \xTo2{f_{lk}} & c_{lk}}\right\},\eqlabel{nervediag}$$
where $c_{ij}\in \C$ are objects, the horizontal arrows $f_{ij}$ are morphisms in $\C$, and the vertical arrows are isomorphisms in~$\C$.
These diagrams stack horizontally, from which one deduces that the resulting simplicial space satisfies the Segal condition.
Furthermore, the space we get from setting $l=0$ consists of chains of invertible morphisms; unraveling the definitions (and with a bit of work), this verifies that the Segal space is complete.

To generalize this construction for a nerve valued in smooth $\infty$-categories, we consider diagrams like \eqref{nervediag} with each $c_{ij}$ an object in a category-valued presheaf on $\Cart$. 

\begin{defn} 
Given a (strict) presheaf $\C\colon \Cart^\op\to \Cat$ of categories, consider the presheaf of groupoids $\C_\Delta\colon \Cart^\op\times \Delta^\op\to \Grpd$ defined by the formula
$$(S,[k])\mapsto {\rm Iso}(\C(S)^{[k]}).$$
Define the \emph{nerve of a category-valued presheaf} as $\N^\Ci(\C):={\rm N}(\C_\Delta)$. 
\end{defn}

\begin{lem}
If $\C$ satisfies descent, then so does $\C_\Delta$ and therefore $\N^\Ci(\C)$ is fibrant in the model structure of \cref{model.smooth.cat}.
\end{lem}

\begin{proof}
This follows immediately from \cref{fibrant.smooth.cat}.
The part with fixed $S$ is verified by Rezk.
The part with fixed $[k]$ is satisfied because $(-)^{[k]}$ and ${\rm Iso}(-)$ both preserve homotopy limits of categories and therefore
preserve the descent property for presheaves of categories.
\end{proof}

We will also need a symmetric monoidal version of the previous lemma.
This is the classical construction of a $\Gamma$-object from a symmetric monoidal object with a strict monoidal structure.

\begin{defn} 
Given a presheaf $\C\colon \Cart^\op\to \SymCat$ of symmetric monoidal categories with a strict monoidal structure (and symmetric strict monoidal functors as morphisms),
consider the presheaf $$\C_\Gamma\colon \Cart^\op\times \Gamma^\op\to \Cat$$ defined by the formula
$$\C_\Gamma(S,\[l])=\C(S)^l$$
with the naturality in $S\in \Cart$ induced by~$\C$ and naturality in $\[l]\in \Gamma$ induced by the strict monoidal structure on $\C(S)$ and its symmetric braiding.
Define the \emph{nerve of a symmetric monoidal category-valued presheaf} by $\N^\Citensor(\C):={\rm N}((\C_\Gamma)_\Delta).$
\end{defn}

We say a presheaf of symmetric monoidal categories with strict monoidal structure
on $\Cart$ satisfies descent if its underlying presheaf of categories (forgetting the symmetric monoidal structure) does. 
This construction continues to work when $\Cart$ is replaced by any other site.
In \cref{cinerve}, we need to use the site $\Cart\times\Gamma$
given by the product of the site $\Cart$
and the category $\Gamma$ equipped with the trivial Grothendieck topology.

\begin{lem}
\label{cinerve}
If $\C$ satisfies descent, then so does $\C_\Gamma$ and therefore $\N^\Citensor(\C)$ is fibrant in the model structure of \cref{model.smooth.mon.cat}.
\end{lem}

\begin{proof}
The functor $(-)^l$ preserves homotopy limits and therefore preserves the descent condition.
Thus $\C_\Gamma$ satisfies descent and we can invoke the previous lemma.
\end{proof}

\section{Bordisms, path categories, and field theories}
\label{sec:1bord}

In this section we present our definition of a smooth 1-dimensional topological field theory over a manifold, as well as a closely related notion of a transport functor.
In \cref{sec:smvect}, we apply the nerve construction of the previous section to the stack of vector bundles on $\Cart$ to obtain the symmetric monoidal smooth $\infty$-category of smooth vector spaces, which is a fibrant object in $\CCat^\otimes$. 
The 1-dimensional oriented bordism ``category'' over $X$, denoted by $1\Bord^\fr(X)$, is defined in \cref{sec:borddef} and is a \emph{non}-fibrant object of $\CCat^\otimes$.
Together, these define the main ingredients of a smooth field theory. 

The model categories $\CCat$ and $\CCat^\otimes$ are simplicial model categories,
so they are equipped with mapping simplicial set functors of the form $$C^\op\times C\to\sSet.$$
These functors are right Quillen bifunctors, and ``derived mapping space'' below refers to their right derived functors.

\begin{defn}
\label{field.theory}
A {\it 1-dimensional smooth oriented topological field theory over~$X$\/} is a point in the derived mapping space $\CCat^\otimes(1\Bord^\fr(X),\Vect)$, and the \emph{space of 1-dimensional smooth oriented topological field theories over~$X$} is this derived mapping space. 
\end{defn}

We will also define a closely related (nonfibrant) smooth $\infty$-category of smooth paths in~$X$, denoted by $\path(X)$, which is a nonfibrant object of $\CCat$. 

\begin{defn}
\label{transport.functor}
A {\it transport functor on~$X$} is a point in the derived mapping space $\CCat(\path(X),\Vect)$, and the {\it space of transport functors on~$X$}
is the derived mapping space $\CCat(\path(X),\Vect)$.
\end{defn}

There is a restriction functor
$$\CCat^\otimes(1\Bord^\fr(X),\Vect)\to \CCat(\path(X),\Vect)$$
from field theories to transport functors. 

Before jumping into detailed definitions of the objects in $\CCat$, we overview some of the ideas that go into
in Lurie's definition of the bordism category~\cite{Lurie:TFT} as a Segal space; see also Calaque--Scheimbauer~\cite{CalaqueScheimbauer}.
The standard way of chopping up a manifold~$M$ into $k$~pieces is a Morse function on~$M$ with a choice of~$k$ nondegenerate critical values.
One can consider the group of diffeomorphisms of~$M$ that preserve the inverse images of these critical values.
The classifying space of this diffeomorphism group is roughly the value of Lurie's Segal space on~$[k]\in \Delta$. 
This is rough because one also wants the classifying space to allow for varying Morse functions and varying critical values. 
In total, the result is a classifying space of the ways of cutting a manifold~$M$ along~$k$ codimension~1 submanifolds. 
In terms of the usual description of the bordism category, this corresponds to a composable chain of bordisms of length $k-1$. 
Forgetting a codimension~1 submanifold or adding extra multiplicity gives maps between these classifying spaces, 
producing the requisite simplicial maps.
In the 1-dimensional case, these classifying spaces are particularly easy owing to the simplicity of Morse decompositions of 1-manifolds. 
Our approach to~$1\Bord^\fr(X)$ is the same idea, but we replace the classifying \emph{space} with a (smooth) classifying \emph{stack}. 

\subsection{Smooth vector spaces}\label{sec:smvect}

Define the sheaf of symmetric monoidal categories, $\V\colon \Cart^\op\to \SymCat$ as follows.
The objects of $\V(S)$ are elements of~$\N$, corresponding to the dimension of a trivial bundle on $S\cong \R^m$, 
and morphisms $m\to n$ are smooth maps $S\to \Hom(\R^m,\R^n)$ into the space of linear maps.
The morphism of groupoids associated with maps $S\to S'$ of objects in $\Cart$ is the identity map on objects, and on morphisms we precompose.
The strict monoidal structure is determined by multiplication in $\N$ and tensor products of linear maps.
The (nontrivial) braiding is induced by the obvious block matrix.
This satisfies descent because smooth functions do. 

\begin{defn}
\label{vector.space.category}
Let $\Vect\in \CCat^\otimes$ (\cref{model.smooth.mon.cat}) be the fibrant object obtained by applying \cref{cinerve} to the sheaf of symmetric monoidal categories $\V$ defined above.
\end{defn}

To explain this a bit more concretely, the vertices of the simplicial set associated to $S\in \Cart$, $[k]\in \Delta$ and $\langle 1\rangle\in \Gamma$ 
are chains of length $k$ of morphisms of vector bundles over~$S$
$$\{V_0\buildrel\phi_1\over\To V_1\buildrel\phi_2\over\To \cdots \buildrel\phi_k\over\To V_k\mid V_i\to S\}\eqlabel{vectob}$$
where the dimensions of the~$V_i$ correspond to the natural numbers in the formal definition.
The 1-simplices of this simplicial set are commutative diagrams of vector bundles:
$$\left\{\cd{V_0&\buildrel\phi_1\over\To& V_1&\buildrel\phi_2\over\To& \cdots &\buildrel\phi_k\over\To& V_k \cr
\mapdown{}&&\mapdown{}&&&&\mapdown{}\cr
V_0'&\buildrel\phi_1'\over\To& V_1'&\buildrel\phi_2'\over\To& \cdots &\buildrel\phi_k'\over\To& V_k'\cr}\right\}\eqlabel{vectmor}$$
where the vertical arrows are vector bundle isomorphisms. 
These 0- and 1-simplices vary with $[k]\in \Delta$ by composing horizontal morphisms of vector bundles or by inserting a horizontal identity morphism of vector bundles.
We can also pull this data back along smooth maps $S'\to S$.

\subsection{The definition of the 1-dimensional bordism category}\label{sec:borddef}

\begin{figure}
\label{fig:bord}
\centerline{\vbox{\newbox\picbox\setbox\picbox\hbox{\includegraphics{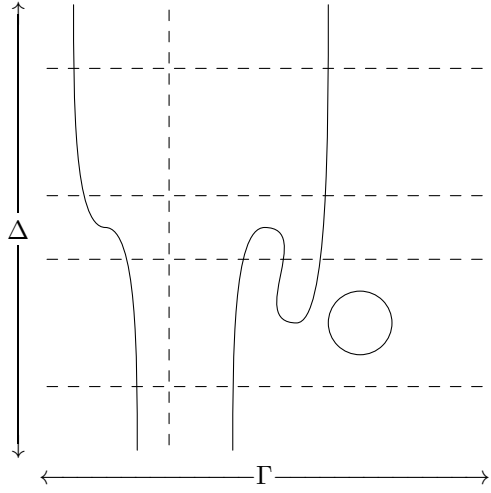}}
\def\sm#1{\vbox{\hbox to0pt{\hss\kern.4pt $#1$\hss}}}
\hbox{\vbox to\ht\picbox{\sm\uparrow\kern-1ex\leaders\vrule\vfil\sm\Delta\kern.5ex\leaders\vrule\vfil\kern-2ex\sm\downarrow}
\kern1ex\vtop{\copy\picbox\hbox to\wd\picbox{\leftarrowfill$\Gamma$\rightarrowfill}}}}}
\caption{A vertex (i.e., an object) in the simplicial set $1\Bord_\pt(\pt)(2,3)$.
The bordism is drawn in solid black, the height function is given by the height in the picture, and the cut functions are represented by the dotted horizontal lines.
The map to $\{0,1,2\} \in \Gamma$ maps the left of the vertical dotted line to~1 and the right of the vertical line to~2.
The fiber over zero is empty.
Regularity at the cut values means that the intersections of the bordism with the dotted lines are transverse.
Restricting attention to the bordism confined within an adjacent pair of horizontal dotted lines gives the three Segal $\Delta$-maps
and similar restrictions corresponding to the vertical dotted line gives the two Segal $\Gamma$-maps.
The action by~$\Sigma_2$ interchanges the bordisms on the left and right sides of the vertical dotted line.
}
\end{figure}

\begin{defn}[The 1-dimensional oriented bordism category over~$X$]
\label{bordism.category}
Given $X\in\Fun(\Cart^\op,\Set)$ (the most important example of $X$ being the presheaf induced by a smooth manifold),
the nonfibrant smooth symmetric monoidal $\infty$-category $1\Bord^\fr(X)$ (\cref{model.smooth.mon.cat})
is defined by taking the nerve of the following presheaf of groupoids on $\Delta\times\Cart\times\Gamma$:
send $[l]\in \Gamma$, $[k]\in \Delta$, and $S\in\Cart$
to the groupoid whose objects are given by data
\begin{enumerate}
\item $M^1$, an oriented 1-manifold that defines a trivial bundle $M^1\times S\to S$,
\item a map $\gamma\colon M^1\times S\to X\times  \{*,1,\ldots, l\}$,
\item cut functions $t_0\le t_1\le t_2\le \cdots\le t_k\in\Ci(S)$,
\item a proper map $h \colon M^1\times S\to \R\times S$ over $S$ called the \emph{height function} such that for each~$s\in S$ the fiber of $h$ over $s\in S$ has $t_i(s)$ as a regular value for all~$i$.
\end{enumerate}
Isomorphisms in the groupoid are certain equivalence classes defined as follows. 
First consider the set of smooth functions $\epsilon_0,\epsilon_k,\epsilon'_0,\epsilon'_k\colon S\to(0,\infty)$
together with orientation-preserving diffeomorphisms
$$\phi\colon h^{-1}((t_0-\epsilon_0,t_k+\epsilon_k)\times S) \to (h')^{-1}((t'_0-\epsilon'_0,t'_k+\epsilon'_k)\times S)$$
over~$S\times X$ such that $\phi$ restricts to a fiberwise diffeomorphism over~$S\times X$ of the form
$$h^{-1}([t_i,t_j]\times S)\to (h')^{-1}([t'_i,t'_j]\times S)$$ for any $0\le i\le j\le k$.
Two such elements are equivalent if their restrictions
to $h^{-1}([t_0,t_k]\times S)$ coincide.
The quotient sets admit a well-defined composition operation, given by pointwise composition.

Functoriality in~$S$ is given by the composition of $\gamma$, $t$, $h$, and $\phi$ with the given map $S'\to S$.
Functoriality in~$\Gamma$ is given by postcomposing~$\gamma$ with the given map of finite sets $\{*,1,\ldots,l\}\to\{*,1,\ldots,l'\}$.
Functoriality in~$\Delta$ with respect to a morphism of simplices $[k']\to[k]$ is given by dropping those~$t_i$ for which $i$ is not in the image of~$[k']$
and duplicating those $t_i$ that are in the image of more than one element of~$[k']$, and restricting~$\phi$ accordingly.
\end{defn}

We observe that $1\Bord^\fr(X)$ is covariant in~$X$: a smooth map $X\to Y$ induces a smooth symmetric monoidal functor $1\Bord^\fr(X)\to 1\Bord^\fr(Y)$.

\begin{rmk}
The properness assumption on 0-simplices in property~(1) guarantees that the 1-dimensional bordism ``between'' $S\times \{t_0\}$ and $S\times \{t_k\}$
is compact in each fiber over~$S$.
\end{rmk} 

\begin{rmk} 
The second piece of data, $M^1\times S\to X\times \{*,1,\dots, l\}$, encodes both the map from the bordism to~$X$, and the partition of connected components of this bordism associated with the monoidal structure.
Our definition of $1\Bord^\fr(X)$ does not satisfy the Segal $\Gamma$-condition, but this failure of fibrancy is not a problem for computing field theories. 
\end{rmk}

\begin{rmk}
\label{modified.bordism}
Below, we will find it convenient to replace $1\Bord^\fr(X)$ with a weakly equivalent object $1\Bord'^\fr(X)$,
which coincides with $1\Bord^\fr(X)$ for all $[n]\in\Delta$
except for $n=0$, where we replace the $\Gamma$-object $1\Bord^\fr(X)(0)$ with the $\Gamma$-object $\[m]\mapsto (X\sqcup X)^m$, where $\[m]=\{*,1,\ldots,m\}$
and $X^m$ denotes the representable presheaf of the $m$-fold cartesian product of~$X\sqcup X$, corresponding to the two orientations of points.
There are canonical homotopy equivalences $1\Bord'^\fr(X)\to1\Bord^\fr(X)$ and $1\Bord^\fr(X)\to1\Bord'^\fr(X)$
that identify points in $X$ with constant paths in~$X$ equipped with the cut function~$t_0=0$.
The advantage of this replacement is that $1\Bord'^\fr(X)(0)(m)$ is a representable presheaf for all $\[m]\in\Gamma$,
hence a cofibrant object in the projective model structure on $\Stacks$.
\end{rmk}

\subsection{The category of smooth paths in a smooth manifold}
Similar to the intuition behind cutting bordisms, we can also consider an analogous structure for paths in~$X$.
In this case a Morse function is afforded by the parametrization of the path itself. 

\begin{defn}
Given $X\in\Fun(\Cart^\op,\Set)$ (the most important example of $X$ being the presheaf induced by a smooth manifold),
define the \emph{smooth path category $\path X\in \CCat$ (\cref{model.smooth.cat}) of~$X$} as the nerve of 
the presheaf of groupoids on $\Delta\times\Cart$ constructed as follows.
A pair $[k]\in \Delta$, $S\in \Cart$
is sent to the groupoid whose objects
consist of a map $\gamma \colon S\times \R\to X$ and cut functions $t_0\le t_1\le \cdots \le t_k \in\Ci(S)$.

Isomorphisms in the groupoid are equivalence classes of a certain equivalence relation.
Elements in the underlying set of this equivalence relation are smooth functions $\epsilon_0,\epsilon_k,\epsilon'_0,\epsilon'_k\colon S\to(0,\infty)$
together with orientation-preserving diffeomorphisms
$$\phi\colon (t_0-\epsilon_0,t_k+\epsilon_k)\times S \to (t'_0-\epsilon'_0,t'_k+\epsilon'_k)\times S$$
over~$S\times X$ such that $\phi$ restricts to a fiberwise diffeomorphism over~$S\times X$ of the form
$$[t_i,t_j]\times S\to [t'_i,t'_j]\times S$$ for any $0\le i\le j\le k$.
Two such elements are equivalent if their restrictions
to $[t_0,t_k]\times S$ coincide.
The quotient sets admit a well-defined composition operation, given by pointwise composition.

Functoriality in~$S$ is given by the appropriate composition of $\gamma$, $t$, and $\phi$ with the given map $S'\to S$.
Functoriality in~$\Delta$ with respect to a morphism of simplices $[k']\to[k]$ is given by dropping those~$t_i$ for which $i$ is not in the image of~$[k']$
and duplicating those $t_i$ that are in the image of more than one element of~$[k']$, and restricting~$\phi$ accordingly.
\end{defn}

We recall that a manifold $X$ defines a presheaf (of sets) on $\Cart$,
and the fiber of~$\path X$ over~$[0]\in\Delta$ is homotopy equivalent to this presheaf via the map~$(\gamma,t_0)\mapsto\gamma(t_0)$,
with an inverse that sends a map $f\colon S\to X$ to an $S$-family of constant paths, $S\times \R\to S\stackrel{f}{\to}X$, with cut function the zero function.
The fiber of $\path X$ over $S\in \Cart$ and $[k]\in \Delta$ is the nerve of the groupoid of $S$-families of paths in~$X$
with $k+1$ marked points and diffeomorphisms of these paths.
The object~$\path X$ is covariant in~$X$, meaning a smooth map $X\to Y$ induces a smooth functor $\path X\to \path Y$,
hence $\path$~is a functor from~$\Mfld$ to~$\CCat$.

There is a smooth functor $\path X\to U(1\Bord^\fr(X))$ in $\CCat$
(where $U$ denotes the forgetful functor of \cref{forget.mon})
we get by viewing a family of paths as the family of bordisms~$S\times M^1=S\times \R$ and the height function~$h$ the projection to~$\R$.
This has an induced restriction map $$\CCat^\otimes(1\Bord^\fr(X),\Vect)\to \CCat(\path X, \Vect)$$
from 1-dimensional field theories over $X$ to representations to the smooth path category of~$X$.
Here $\CCat^\otimes(-,-)$ and $\CCat(-,-)$ denote the corresponding {\it derived\/} mapping simplicial sets.

\begin{rmk}
\label{modified.path}
In analogy to \cref{modified.bordism},
we will find it convenient to replace $\path X$ with a weakly equivalent object $\path' X$,
which coincides with $\path X$ for all $[n]\in\Delta$
except for $n=0$, where we replace $\path X(0)$ with $X$ itself.
There are canonical homotopy equivalences $X\to\path X(0)$ and $\path X(0)\to X$ that identify $X$ with constant paths in~$X$ equipped with the cut function~$t_0=0$.
The advantage of this replacement is that $\path' X(0)$ is a representable presheaf, hence a cofibrant object in the projective model structure on $\Stacks$.
\end{rmk}

\subsection{Descent for field theories and representations of paths}

A key step to verifying the main theorem is that field theories over~$X$ and representations of smooth paths in~$X$ can be computed locally in the following sense. 

{\def\thethm{C}
\begin{thm}
\label{thm:local}
Let $\{U_i\}$ be an open cover of a smooth manifold~$X$.
The canonical maps
$$\hocolim \path(U_k)\to \path(X),\qquad \hocolim 1\Bord^\fr(U_k)\to 1\Bord^\fr(X)$$
are equivalences of smooth $\infty$-categories and symmetric monoidal smooth $\infty$-categories, respectively.
Here $k$ runs over all finite tuples of elements in~$I$ and $U_k$ denotes the intersection of $U_i$ for all $i\in k$.
This immediately implies that the assignments $$X\mapsto \CCat(\path(X),\Vect),\qquad X\mapsto \CCat^\otimes(1\Bord^\fr(X),\Vect)$$ are stacks on the site of smooth manifolds. 
\end{thm}
}

\begin{proof}
See Grady--Pavlov \cite[Theorem~1.0.1]{GradyPavlov.Local} for the case of bordism categories.
The case of path categories then follows formally.
We remark that the most technical part of the cited proof (Section~6.6 in op.~cit.)\ becomes completely trivial in the 1-dimensional case,
since the nerves of relevant categories are contractible for trivial reasons.
\end{proof}

We apply the above result to reduce our main theorems to the case $X\in\Cart$.
In particular, this simplifies the construction of a transport functor from a vector bundle with connection,
since the general case $X\in\Mfld$ would require us to work with arbitrary cocycles for vector bundles,
bringing considerable technicalities,
whereas for $X\in\Cart$ all vector bundles over~$X$ are trivial, and the problem reduces to manipulating connection 1-forms.

\begin{defn}
\label{vect.connection}
We define $\Vect^\nabla\in\Stacks$ as follows.
Given $X\in\Cart$, we send it to the nerve of groupoid
whose objects are pairs $(n,\omega)$,
where $n\ge0$ specifies a finite-dimensional vector space $V=\R^n$ and $\omega\in\Omega^1(X,\End(V))$.
Morphisms $(n,\omega_0)\to(n,\omega_1)$ are smooth maps $f\in \Ci(X,\GL(V))$
such that $\omega_1=\Ad_{f^{-1}}\omega_0+f^{-1} df$,
where $\Ad$ denotes the adjoint action.
\end{defn}

The groupoid $\Vect^\nabla(X)$ is equivalent to the groupoid of trivial vector bundles with connection over the cartesian space~$X$.

\begin{cor}
\label{connection.transport}
Consider the functors
$$\CCat(-,\Vect^\nabla), \CCat(\path (-), \Vect)\colon \Fun(\Cart^\op,\Set)^\op \to \sSet$$
that send $X\in\Fun(\Cart^\op,\Set)$ to $\CCat(X,\Vect^\nabla)$ and $\CCat(\path X, \Vect)$, respectively.
The space of natural weak equivalences
$$\CCat(-,\Vect^\nabla) \to \CCat(\path (-), \Vect)$$
is weakly equivalent to the space of natural weak equivalences
between the same functors, restricted along the Yoneda embedding $\Cart\to\Fun(\Cart^\op,\Set)$ to the category~$\Cart$.
In particular, any natural weak equivalence
$$\Vect^\nabla \to \CCat(\path (-), \Vect)$$
on $\Cart$ can be extended to a natural weak equivalence on $\Fun(\Cart^\op,\Set)$ in a unique way up to a contractible choice.
\end{cor}

\begin{proof}
$\Cart$ generates $\Fun(\Cart^\op,\Set)$ under homotopy colimits.
Both functors, $$\CCat(-,\Vect^\nabla) \qquad {\rm and}\qquad \CCat(\path (-), \Vect),$$ send homotopy colimits in $\Fun(\Cart^\op,\Set)$ to homotopy limits in $\sSet$.
For the former, it boils down to the classical fact that vector bundles with connection satisfy descent,
whereas for the latter it follows from \cref{thm:local}.
\end{proof}

In the next section, we construct a specific weak equivalence
$$\Vect^\nabla \to \CCat(\path (-), \Vect)$$
in $\PreStacks$,
which will prove that representations of the path category of any smooth manifold~$X$
are precisely vector bundles with connection over~$X$.

\section{Representations of the smooth path category of a manifold}
\label{sec:paths}

In this section we prove the following. 
\begin{prop}
\label{prop:pathvb}
Given $X\in\Fun(\Cart^\op,\Set)$,
the derived mapping space $$\CCat(\path X, \Vect)$$
is naturally weakly equivalent to $\Vect^\nabla(X)$ of \cref{vect.connection}.
(In particular, we can take $X\in\Mfld$.)
\end{prop}

\begin{proof}
By \cref{connection.transport},
it suffices to construct such a natural weak equivalence for $X\in\Cart$.
Replace $\path X$ with weakly equivalent $\path' X$ from \cref{modified.path}.
Recall that $\Vect$ is fibrant in $\CCat$ (\cref{fibrant.smooth.cat}).
Furthermore,
the stack of objects $(\path' X)_{[0]}=\path' X(0)=X$ is a representable (hence cofibrant) presheaf in $\PreStacks$
and $\Vect$ and $\path' X$ are constructed as objectwise nerves of groupoids.
Hence, the derived mapping space $\CCat(\path' X, \Vect)$ can be computed using the nonderived hom $\CCat(\path' X,\Vect)$.
The map $$\Vect^\nabla(X)\to \CCat(\path' X,\Vect)$$ is constructed in \cref{cartesian.connection.transport} and is shown to be an isomorphism in
\cref{transport.functor.data} and \cref{transport.connection.isomorphism}.
\end{proof}

The following construction codifies the parallel transport data of a connection on a vector bundle as a functor.

\begin{defn}
\label{cartesian.connection.transport}
Given $X\in\Cart$, we construct a map (natural in $X$)
$$\Vect^\nabla(X)\to \CCat(\path' X,\Vect)$$ as follows.
A (trivial) vector bundle with connection~$(n,\omega)$ defines a smooth functor $R\colon \path' X\to \Vect$ via parallel transport: to $f\colon S\to X$
an $S$-family of points in~$X$ we assign the object~$n$ over~$S$, defining a functor $X=\path' X(0)\to \Vect(0)$.
To a family of oriented paths $S\times\R\to X$, we apply the fiberwise parallel transport with respect to the connection 1-form~$\omega$,
yielding a morphism of (trivial) vector bundles over~$S$.
These maps are invariant under families of diffeomorphisms of 1-manifolds, so we obtain a functor $\path'_S X(1)\to \Vect_S(1)$,
which is again natural in~$S$ so defines a fibered functor $\path' X(1)\to \Vect(1)$
by extension from individual fibers in the usual fashion.
We extend in the obvious way to $\path' X(k)\to \Vect(k)$, where naturality with respect to
maps in~$\Delta$ follows from compatibility of parallel transport with concatenation of paths.
Hence, we have constructed a functor from the path category of~$X$ to smooth vector spaces.
Lastly, we observe that an isomorphism of vector bundles with connection leads to a natural isomorphism of functors of such functors, i.e., an edge in the simplicial mapping space.
An $n$-simplex comes from a composable $n$-tuple of isomorphisms of vector bundles with connection.
\end{defn}

\subsection{Reduction to parallel transport data}

In this section we whittle the proof of \cref{prop:pathvb} down to a statement about parallel transport data, by which we shall mean smooth endomorphism-valued
functions on paths that compose under concatenation of paths and are compatible with restrictions to intersections of the cover.
The next definition and lemma describe the precise manner in which a representation of the path category determines a transport functor on~$X$.

\begin{defn}
\label{tran}
We define the stack $\Tran\in\PreStacks$ of transport data as follows.
Given $X\in\Cart$, we send it to the nerve of groupoid whose objects are
pairs $(n,F)$, where $n\ge0$ determines a vector space $V=\R^n$ and $F$ is a
morphism $$F\colon\R\times\Hom(\R,X)\to\End(V)$$ in $\PreStacks$
such that the following three properties are satisfied:
(1)~for any $p\colon\R\to X$ and any $L_1,L_2\in \R$ the following functoriality property holds:
$$F(L_2,p\circ S_{L_1})\circ F(L_1,p)=F(L_1+L_2,p),$$
where $S_{L_1}(t)=t-L_1$;
(2)~for any $p\colon\R\to X$ we have $F(0,p)=\id_V$;
(3)~$F$ is invariant under diffeomorphisms of paths:
if $g\colon \R\to\R$ is an orientation-preserving diffeomorphism such that $g(0)=0$, then $$F(L,p)=F(g^{-1}(L),p\circ g).$$
(The first argument of~$F$ specifies the length~$L$ of a smooth path~$p$ in~$X$.
The path itself is given by the second argument of~$F$ and is parametrized by $[0,L]\subset\R$.)
Morphisms $F_1\to F_2$ are smooth maps $h\colon X\to\GL(V)$ such that $$h(p(L))\circ F_1(L,p)=F_2(L,p)\circ h(p(0)).$$
\end{defn}

\begin{lem}
\label{transport.functor.data}
There is an isomorphism in $\PreStacks$
$$\CCat(\path' (-),\Vect)\to \Tran.$$
($\path' X$ is constructed in \cref{modified.path}.)
\end{lem}

\begin{proof}
Fix some $X\in\Cart$; we need to define a morphism $$\CCat(\path' X,\Vect)\to\Tran(X).$$
Both sides are nerves of groupoids, so we define the map first on objects, then on morphisms.
Pick a functor $R\colon\path' X\to\Vect$.
The presheaf of objects $\path X(0)=X$ maps via~$R$ to a fixed object of $\Vect$ given by some dimension~$n\ge0$.
The data of~$F$ is obtained by evaluating on $[1]\in\Delta$ with cut function $t_0=0$.
On objects, we map $R\mapsto(n,F)$.
Property~(1) boils down to functoriality with respect to the three coface maps $[1]\to[2]$ in $\Delta$,
where the middle face map computes the composition.
Property~(2) boils down to functoriality with respect to the codegeneracy map $[1]\to[0]$ in $\Delta$.
Property~(3) boils down to the fact that isomorphisms in $\path X([n],S)$ are endpoint-preserving diffeomorphisms between $n$-chains of paths in~$X$,
whose individual $(n+1)$ vertices are identity maps (this holds for $\path' X$, not for $\path X$).
Finally, a morphism $R\to R'$ is given by a map $h\colon X=\path' X(0)\to \GL(n)$,
which yields a morphism $h\colon(n,F)\to(n',F')$ (where $n=n'$).

Conversely, the inverse map $$\Tran(X)\to \CCat(\path' X,\Vect)$$
sends $(n,F)$ to the functor $R\colon\path' X\to \Vect$
that maps the presheaf of objects $\path' X(0)$ to~$n$.
The presheaf of $k$-simplices for $k\ge1$
is mapped to the corresponding transport maps between endpoints.
\end{proof}

\subsection{From parallel transport data to vector bundles with connection}

From the above discussion, we have shown that a point in the (derived) mapping space~$\CCat(\path X, \Vect)$ for $X\in\Cart$
defines a parallel transport data for a vector bundle on~$X$, i.e., an object of $\Tran(X)$.
In this section we explain how parallel transport data defines a vector bundle with connection.
More precisely, we construct an equivalence $\Tran(X)\to\Vect^\nabla(X)$.
Most of the ideas below are present in Freed~\cite{Freed} and Schreiber--Waldorf~\cite{SchreiberWaldorf},
and we have adapted them to our situation with some minor modifications.

\begin{lem}
Given $X\in\Cart$ and $(n,F)\in\Tran(X)$,
the map~$F$ assigns the identity map on $V=\R^n$ to constant paths in~$X$.
\end{lem}

\begin{proof}
Since a constant path~$\gamma$ can be factored as the concatenation~$\gamma*\gamma$,
the value of~$F$ on~$\gamma$ must be a projection in~$V$, denoted by~$P_\gamma$.
Furthermore, there is a family of constant paths parametrized by $[0,t]$ coming from the restriction of~$\gamma$ to $[0,t']\subset [0,t]$.
Over~$t'=0$, the constant path is the identity morphism in the path category
and therefore is assigned the identity linear map.
Smoothness gives a family of projections connecting~$P_\gamma$ on~$V$ that is the identity projection at an endpoint.
Since the rank of the projection is discrete, it must be constant along this family.
Therefore, $P_\gamma$~is the identity.
\end{proof}

\begin{lem}
Given $X\in\Cart$ and $(n,F)\in\Tran(X)$,
the map~$F$ lands in the invertible morphisms, i.e.,
the morphism $F(\gamma)$ for an family of paths~$\gamma$ is an \emph{isomorphism} on $V=\R^n$.
\end{lem}

\begin{proof}
Since a path of length zero is assigned the identity linear map on~$V$, by continuity there is an $\epsilon>0$ such that
the restriction of any path~$\gamma$ to $[0,\epsilon]$ is assigned an invertible morphism.
Observe that this holds for any point on a given path (though possibly with variable~$\epsilon$).
Choosing a finite subcover and
factoring
the value on a path into the value on pieces of the path subordinate to the subcover,
we see that the value on a path is a composition of vector space isomorphisms, and therefore an isomorphism.
\end{proof}

The remaining work is in the construction of an inverse to the map $\Vect^\nabla(X)\to\Tran(X)$
(see \cref{cartesian.connection.transport} and \cref{transport.functor.data}).
For this we use the following two lemmas of Schreiber--Waldorf~\cite[Lemmas 4.1 and~4.2]{SchreiberWaldorf}, reproduced here for convenience.

\begin{lem}
\label{KonradUrs1}
For a finite-dimensional vector space~$V$,
smooth functions $$F\colon \R\times \R\to\Aut(V)$$
satisfying the cocycle condition $F(y,z)\cdot F(x,y)=F(x,z)$ and $F(x,x)=\id$ are in bijection with 1-forms, $\Omega^1(\R;\End(V))$.
\end{lem}

\begin{proof}
Given such a 1-form~$A$, consider the initial value problem
$$(\partial_t\alpha)(t)=A_t(\partial_t)(\alpha(t)),\qquad \alpha(s)=\id,\eqlabel{eq:IVP}$$
where $\alpha\colon \R\to \Aut(V)$ and $s\in \R$.
We obtain a unique solution~$\alpha(t)$ depending on~$s$, and define~$F(s,t)=\alpha(t)$.
The function~$F$ is smooth in~$s$ because the original coefficients were smooth in~$s$
and is globally defined because the equation is linear.
To verify that $F(s,t)$ satisfies the cocycle condition, we calculate
$$\partial_t \left(F(y,t)F(x,y)\right)=\left(\partial_t F(y,t)\right) F(x,y)=A_t(\partial_t)F(y,t)F(x,y),$$
and since $F(y,y)F(x,y)=F(x,y)$, uniqueness dictates that $F(y,t)F(x,y)=F(x,t).$
Conversely for $F\colon \R\times \R\to \Aut(V)$,
let $\alpha(t)=F(s,t)$ for some $s\in \R$ and let
$$A_t(\partial_t)=\left(\partial_t\alpha(t)\right)\alpha(t)^{-1}.$$
When $F$ satisfies the cocycle condition, $A_t(\partial_t)$ is independent of the choice of~$s$:
$$F(s_0,t)=F(s_1,t)F(s_0,s_1)\ \implies \ (\partial_t F(s_0,t) )F(s_0,t)^{-1}=(\partial_t F(s_1,t))F(s_1,t)^{-1}.$$
This gives the desired bijection.
\end{proof}

\begin{lem}
\label{KonradUrs2}
Let $A,A'\in \Omega^1(\R;\End(V))$ be endomorphism valued 1-forms on~$\R$ and $g\colon \R\to \Aut(V)$
be a smooth function. Let~$F_A$ and~$F_{A'}$ be the smooth functions corresponding to~$A$ and~$A'$ by \cref{KonradUrs1}. Then 
$$g(y)\cdot F_A(x,y)=F_{A'}(x,y)\cdot g(x)$$
if and only if $A'=\Ad_{g^{-1}} A+g^{-1}dg$.
\end{lem}

\begin{proof}
The function $g(y) F_A(x,y) g(x)^{-1}$ solves the initial value problem~\eqref{eq:IVP} for~$A'$:
$$\eqalign{\partial_y(g(y)F(x,y)g(x)^{-1})&=(\partial_yg(y))F(x,y)g(x)^{-1}+g(y)\partial_yF(x,y)g(x)^{-1}\cr
&=(\partial_yg(y)g(y)^{-1}) (g(y)F(x,y)g(x)^{-1})\cr
&\qquad{}+ (g(y)A_y(\partial_y)g(y)^{-1})(g(y)F(x,y)g(x)^{-1})\cr}$$
so by uniqueness we obtain the desired equality.
\end{proof}

Now we construct a differential form from the parallel transport data that will give rise to a connection.
Throughout, $(d,F)\in\Tran(X)$ is a transport data on~$X\in\Cart$ with typical fiber $V=\R^d$.
Let $\gamma\colon \R\to X$ be a path such that $\gamma(0)=p$ and $\dot\gamma(0)=v$;
restrictions of~$\gamma$ to intervals (as a family over~$\R^2$) will give a family of paths in~$\path X$, i.e., a 0-simplex.
Define $$F_\gamma(x,y)=F(\gamma\colon [x,y]\to X ), \quad F\colon \R\times \R \to\End(V).$$
By the above lemma, $F_\gamma$~gives us a 1-form~$A_\gamma$ with values in~$\End(V)$.
By varying~$\gamma$, we want to promote this to a 1-form on~$X$ whose value at $(p,v)$ is~$A_\gamma(\partial_t)$.

\begin{defn}
\label{disintegrate}
The morphism $\Tran\to\Vect^\nabla$ is defined as follows.
Fix $X\in\Cart$ and $(d,F)\in\Tran(X)$; we need to produce $A\in\Omega^1(X,\End(V))$, where $V=\R^d$.
Fix a point $p\in X$ and a tangent vector $v\in X$; we need to define $A_p(v)\in\End(V)$.
We use the linear structure on~$X\cong\R^n$ to define $A_p(v)=(A_{tv})_0(1)$ where~$tv$ is the path through~$p\in X$ with velocity vector~$v\in X$,
$A_{tv}$ is a 1-form on~$\R$, and $(A_{tv})_0(1)$ is the value of this 1-form at 0 evaluated at~$1\in T_0\R$.
This defines a functor $\Tran(X)\to\Vect^\nabla(X)$ on objects,
and on morphisms we send $h\colon X\to\GL(V)$ to itself, now viewed as a morphism in $\Vect^\nabla(X)$.
\end{defn}

\begin{lem}\label{itsalemma}
\cref{disintegrate} is well-defined: objects in $\Tran$ are sent to smooth differential 1-forms $A\in\Omega^1(X,\End(V))$, and isomorphisms are sent to gauge transformations $A'\mapsto\Ad_{g^{-1}} A+g^{-1}dg$.
\end{lem}

\begin{proof}
The claim on isomorphisms follows from \cref{KonradUrs2}. To verify the claim on objects, we first observe that~$A$ is smooth: choose families of affine paths in a neighborhood of~$p$ and invoke smoothness of the representation.
Furthermore, we claim that~$A$ satisfies~$A(\lambda v)=\lambda A(v)$ for all~$\lambda>0$.
To see this, define $\gamma_\lambda(t)=\gamma(\lambda t)$ for~$\lambda>0$.
We compute $$A(\lambda v)=\partial_t F_{\gamma_\lambda}(0,t)|_{t=0}=\partial_t F_\gamma(0,\lambda t)|_{t=0}=\lambda A(v).$$
\cref{homogeneous.linear} shows that this property implies~$A$ is linear. 
\end{proof}

\begin{lem}
\label{homogeneous.linear}
A smooth function $A\colon V\to W$ between vector spaces that satisfies $A(\lambda v)=\lambda A(v)$ for $\lambda>0$ is linear.\label{lem:lin}
\end{lem}

\begin{proof}
It suffices to show that $A$~is equal to its derivative at zero.
From the assumptions it follows that~$A(0)=0$.
Smoothness of~$A$ implies that $dA(0)$ exists, and we compute its value on~$v$ by the one-sided limit:
$$(dA(0))(v)=\lim_{\lambda\to 0^+} A(\lambda v)/\lambda=\lim_{\lambda\to 0^+} \lambda A(v)/\lambda=A(v),$$
completing the proof.
\end{proof}

The next lemma shows that~$A$ determines the given representation.
Our techniques are in the spirit of D.~Freed's \cite[Appendix~B]{Freed}, though benefited from K.~Waldorf pointing out to us the utility of Hadamard's lemma in this context.

\begin{lem}
\label{integrate}
For $X\in \Cart$ and $(d,F)\in\Tran(X)$,
the value of~$F$ on a path~$\gamma$ is
the path-ordered exponential associated to the $\End(V)$-valued 1-form~$A$ constructed in \cref{disintegrate}, where $V=\R^d$.
\end{lem}

\begin{proof}
Let $X\cong \R^n$ and $\gamma\colon [0,T]\to \R^n$ be a path.
Fix $N$ a large integer, and let $\gamma_i$ denote the restriction of~$\gamma$ to $[T(i-1)/N,Ti/N]$ for $1\le i\le N$.
By definition of $\Tran(X)$,
$$F(\gamma)=F({\gamma_N})\circ \cdots \circ F({\gamma_2})\circ F({\gamma_1}).$$
Reparametrize~$\gamma_i$ by~$\tilde\gamma_i(t):=\gamma_i(T(t+i-1)/N)$ and let~$\ell_i\colon [0,T/N]\to \R^n$ denote the affine path of length~$1$
starting at~$\gamma_i(0)$ with velocity~$\dot{\gamma}_i(0)=v_i$.
By Hadamard's lemma there is a smooth function~$g_i$ with $\gamma_i(s)-\ell_i(s)=s^2 g_i(s)$.
Define $G\colon [0,1]\to \End(V)$ by $G(t):=F(\tilde\gamma_i|_{[0,t]})$.
Using that $\tilde{\gamma}_i(s)= \ell_i(s(T/N))+s^2(T^2/N^2)g_i(s(T^2/N^2))$ and applying Hadamard's lemma to~$G$ we obtain
$$G(t)=G(0)+t G'(0)+t^2G_2(t)=\id+t(T/N) A_{\ell_i}(v_i)+O(N^{-2})$$
for some function $G_2\colon [0,1]\to \End(V)$.
The~$O(N^{-2})$ estimate comes from Taylor's formula and the fact that the original domain of definition~$[0,T]$ is compact,
so a uniform estimate can be given for the coefficient before~$(T/N)^2$.
The claimed form of the derivative~$G'(0)$ follows from \cref{lem:lin} and an argument in Schreiber--Waldorf~\cite[Lemma~B.2]{SchreiberWaldorf}
(reproduced in the next paragraph) to show that $A_{\gamma_i}(v_i)=A_{\ell_i}(v_i)$.

First we consider the family of paths $\Gamma(t,\alpha):=\ell_i(t)+\alpha g_i(t)$ depending on the parameter~$\alpha$ for~$0\le \alpha \le 1$.
Define $q\colon [0,1]^2\to [0,1]^2$ by $(t,\alpha)\mapsto (t,t^2\alpha)$.
The composition $$(\Gamma\circ q)(t,\alpha)=\ell_i(t)+\alpha t^2g_i(t)$$
defines a smooth homotopy running from $\ell_i$ (when $\alpha=0$) to $\gamma_i$ (when $\alpha=1$).
For a fixed $\alpha$, we evaluate $F$ on the family of paths $\Gamma\circ q$ obtained from restriction to $[0,t]\times \{\alpha\}\subset [0,1]^2$
and differentiate with respect to $t$ using the chain rule,
$${d\over dt}F((\Gamma\circ q)|_{[0,t]\times\{\alpha\}})|_{t=0}=d(F(\Gamma))|_{q(0,\alpha)}\circ{dq\over dt}\Big|_{t=0}=d(F(\Gamma))|_{(0,0)}\circ (1,0).$$
The right hand side is independent of $\alpha$, whereas the left hand side is $A_{\ell_i}(v_i)$ for $\alpha=0$ and $A_{\gamma_i}(v_i)$ when~$\alpha=1$, so the claim follows.

Putting this together, we have $F(\gamma_i)=\id+(T/N) A_{\ell_i}(v_i)+O(N^{-2})$,
and taking $N\to\infty$,
$$\eqalign{
F(\gamma)&=\lim_{N\to \infty} (\id+(T/N) A_{\ell_1}(v_1))(\id+(T/N) A_{\ell_2}(v_2))\cdots (\id+(T/N) A_{\ell_N}(v_N))\cr
&=\lim_{N\to \infty} \exp((T/N)A_{\ell_1}(v_1)) \exp((T/N)A_{\ell_2}(v_2))\cdots \exp((T/N)A_{\ell_N}(v_N))\cr
&=\path \exp(A(\dot{\gamma})),\cr
}$$
since the limit is the definition of the path-ordered exponential of~$A$ along~$\gamma$.
\end{proof}

\begin{prop}
\label{transport.connection.isomorphism}
The morphism $$\Tran\to\Vect^\nabla$$ in $\PreStacks$ constructed in \cref{disintegrate} is an isomorphism.
\end{prop}

\begin{proof}
The inverse isomorphism is the composition $$\Vect^\nabla(X)\to\CCat(\path' X,\Vect)\to\Tran(X),$$
where the two maps are constructed in \cref{cartesian.connection.transport} and \cref{transport.functor.data}.
As shown in \cref{integrate}, for every $X\in\Cart$, the composition $$\Tran(X)\to\Vect^\nabla(X)\to\Tran(X)$$ is an isomorphism on objects.
Morphisms in $\Tran(X)$ and $\Vect^\nabla(X)$ were defined as smooth maps $X\to\GL(V)$ satisfying certain respective properties,
and we have shown in \cref{itsalemma} that these properties are preserved by these functors.
\end{proof}

\numberwithin{equation}{section}

\section{Smooth 1-dimensional field theories and the cobordism hypothesis}
\label{sec:1TFT}

We keep the notation of the previous section: $X\in\Cart$ is a cartesian space,
on which we consider a field theory.

By construction, there is a map~$\path X\to U(1\Bord^\fr(X))$
(in the category $\CCat$ of \cref{model.smooth.cat}, where $U$ is the forgetful functor of \cref{forget.mon})
that views a path as a bordism.
This will allow us to apply arguments from the preceding section to the bordism category.

\begin{lem}
\label{lem:sitting}
The value of a field theory (see \cref{field.theory}) on a family of bordisms $\gamma\colon S\times M\to X$ as a vertex in $1\Bord^\fr_S(X)$
is equal to the value on a bordism $\gamma_{\rm sit}\colon S\times M\to X$ that has the same image in~$X$ as~$\gamma$ but has sitting instants,
meaning that the map $\gamma_{\rm sit}\colon M\to X$ is constant near~$t_0$ and~$t_1$.
\end{lem}

\begin{proof}
Using the $\Gamma$-structure, it suffices to prove the lemma for arcs in~$X$, i.e., $S$-families~$\gamma \colon S\times I\to X$ for~$I$ an interval.
Choose~$b\colon\R\to\R$ to be a smooth bump function such that $b|_{(-\infty,1/3]}=0$, $b|_{[2/3,\infty)}=1$, and $b|_{(1/3,2/3)}\subset(0,1)$.
Consider a new $S\times\R$-family of 1-manifolds
that for $t\in\R$ is given by $\gamma_{\rm sit}:=\gamma\circ \Gamma(x,t)$, where $\Gamma(x,t)=tx+(1-t)b(x)$ for $x\in I$.
To this family a field theory assigns a smooth family of linear maps.
We observe that for all $t\in(0,1]$, the fibers in this family are isomorphic as morphisms in the fiber of $1\Bord^\fr(X)$ over $S=\pt\in \Cart$.
Therefore, a field theory assigns the \emph{same} linear maps for all $t\ne 0$.
By smoothness, we obtain the same linear map at $t=0$
and the resulting path has sitting instants around 0~and~1 by construction.
\end{proof}

We are now ready to prove \cref{thmb}.

{\def\thethm{D}
\begin{thm}
\label{thmb.details}
Evaluating at $\[1]\in\Gamma$ and restricting along $\path X\to U(1\Bord^\fr(X))$ yields a weak equivalence of derived mapping simplicial sets, natural in $X\in\Mfld$
(and, more generally, $X\in\Fun(\Cart^\op,\Set)$):
$$\CCat^\otimes(1\Bord^\fr(X),\Vect^\otimes)\to \CCat(\path X,\Vect).$$
Thus, there is an equivalence between 1-dimensional oriented topological field theories over~$X$ valued in~$\Vect^\otimes$
and \ci-functors from the smooth path category of~$X$ to~$\Vect$.
\end{thm}
}

\begin{proof}
Applying \cref{thm:local}, we reduce the problem to the case $X\in\Cart$.
Applying \cref{modified.bordism}, we replace $1\Bord^\fr(X)$ with $1\Bord'^\fr(X)$,
for which we have $1\Bord'^\fr(X)([0],\[1])=X\sqcup X$, corresponding to two possible orientations of points in~$X$.
These two copies of~$X$ a priori map to some $d,d'\in\N=\Vect(X)([0],\[1])$,
which uniquely determines the maps on $1\Bord'^\fr(X)([0],\[m])$ for $m\ne1$.
As shown below, we necessarily have $d=d'$, which corresponds to the dimension of the vector bundle determined by the field theory.

To understand the value of a field theory on morphisms, since the target category $\Vect^\otimes$ is fibrant (\cref{fibrant.smooth.mon.cat}),
in particular, satisfies the Segal $\Delta$-condition,
a functor $1\Bord^\fr(X)\to\Vect$ is determined (up to a contractible choice) by its value over the fiber $[1]\in\Delta$.
Furthermore, since any bordism can be expressed as a disjoint union of connected bordisms, we can restrict attention to $S$-families of connected 1-manifolds in $1\Bord^\fr_S(X)(1)$.

In the case that cut functions satisfy $t_0<t_1$, Morse theory of 1-manifolds cuts a given connected bordism into elementary pieces that are of three types:
(1)~bordism from a point to a point (all points of $M^1\times \{s\}$ are regular values for~$h$),
(2)~bordisms from the empty set to a pair of points (0-handles), and
(3)~bordisms from a pair of points to the empty set (1-handles).
For a given bordism, i.e., 0-simplex of $1\Bord^\fr(X)$, this reduction comes from a choice of (new)
height function that is Morse with regular values at the prescribed cut values, which defines a 1-simplex in $1\Bord^\fr(X)$
connecting the original bordism to one with a Morse height function.
Then we can impose additional cut points using the Morse height function to reduce to the cases above.
The relations among these generators are precisely the familiar birth-death diagrams from 1-dimensional Morse theory.

When cut functions satisfy $t_0=t_1$, since $t_0$ is a regular value and the bordism is connected, this bordism is in the image of the degeneracy map,
i.e., is an identity path in the bordism category.
For $S$ connected and $t_0\le t_1$ with $t_0=t_1$ somewhere on $S$, then this is necessarily a bordism of type~(1) above.

In the case that the above types of generating bordisms are mapped constantly to~$X$, meaning the map $x\colon M^1\times S\to X$ factors through
the projection to~$S$, the standard dualizable object argument
shows that the value of the field theory on the $(+)$-point must be a vector space~$(V_+)_x=\R^d$, and the value on the $(-)$-point is the dual space,~$(V_-)_x^*=\R^{d'}$,
which in our formulation amounts to showing $d=d'$.

Now we need to show that the value on a generating bordism with an arbitrary map to~$X$
is determined by
the value of the field theory on the path category.
For generating bordisms of type~(1) this is clear, since such a bordism can be identified with a morphism in the path category.

For bordisms of type (2) and~(3) we use \cref{lem:sitting} to identify the value of a field theory on a 0- or 1-handle
with the value on a handle that has a sitting instant at its Morse critical point.
Then we can factor the handle into~3 pieces:
one given by a subset of the sitting instant of the Morse critical points (i.e., a handle that is mapped constantly to~$X$) and two paths given
by the closure of the complement of this subset in the original handle.
Hence, the value of the original bordism is determined by previously computed dualizing data at the sitting instant
together with the value on paths between points.
\end{proof}

\begin{proof}[Proof of \cref{thma}]
The result follows from \cref{thmb} and \cref{prop:pathvb}.
\end{proof}

\section{Reedy model structures}\label{Reedy}

In this auxiliary section we review the necessary facts from the theory of Reedy model structures.

Let $C$ be a model category.
Following Hirschhorn~\cite[Definition~15.3.3]{Hirschhorn}, we review the Reedy model structure on the category of simplicial objects in $C$, i.e., the category of functors $\Delta^\op\to C$. 

First, define the \emph{$n$th latching functor}, $L_n\colon\Fun(\Delta^\op,C)\to C$  as 
$$L_nX=\colim_{[m]\to[n]}X_m,$$ 
where the colimit is indexed by surjections of finite ordered sets $[m]\gets[n]$ that are not isomorphisms (i.e., the union of degenerate simplices).
Similarly, define the \emph{$n$th matching functor} $M_n\colon\Fun(\Delta^\op,C)\to C$ as
$$M_nX=\lim_{[m]\gets[n]}X_m,$$ 
where the limit is indexed by injections of finite ordered sets $[m]\to[n]$ that are not isomorphisms (i.e., the defining data of a boundary of an $n$-simplex). 

Now, in the Reedy model structure on simplicial objects in~$C$, a map $X\to Y$ is a cofibration if 
$$X_n\sqcup_{L_nX}L_nY\to Y_n$$ 
is a cofibration in~$C$ for any $[n]\in\Delta$.
Similarly, a map $X\to Y$ is a fibration if 
$$X_n\to Y_n\times_{M_nY}M_nX$$ 
is a fibration for any $[n]\in\Delta$. 
In particular, an object~$X$ is cofibrant if the latching map $L_nX\to X$ is a cofibration in $C$ for any~$[n]$
and fibrant if the matching map $X\to M_nX$ is a fibration in $C$ for any $[n]\in\Delta$.

By Hirschhorn~\cite[Theorem~15.6.27]{Hirschhorn} the Reedy model is cofibrantly generated,
with generating (acyclic) cofibrations as in Hirschhorn~\cite[Definition~15.6.23]{Hirschhorn}:
if $A\to B$ is a generating (acyclic) cofibration in $C$,
then $A\otimes[n]\sqcup_{A\otimes\partial[n]}B\otimes\partial[n]\to B\otimes[n]$
is a generating (acyclic) cofibration of the Reedy model structure.

%$$

%$$

\end{document}